\documentclass[12pt]{amsart}  
\usepackage{setspace}  
\usepackage{float}
\usepackage{xcolor}
\usepackage{amsmath,amssymb, amsthm,graphicx}
\usepackage{verbatim,enumerate}
\usepackage{subcaption}
\usepackage{mwe}

\usepackage[capitalise]{cleveref} 

\newtheorem{theorem}{Theorem}
\newtheorem{corollary}[theorem]{Corollary}
\newtheorem{lemma}[theorem]{Lemma}

\theoremstyle{remark}
\newtheorem{remark}[theorem]{Remark}

\theoremstyle{definition}

\usepackage{mathtools}

\DeclareMathOperator{\conv}{conv} 

\usepackage{anysize}
\marginsize{3cm}{3cm}{3cm}{3cm}

\begin{document}
\def\red{\textcolor{red}}
\def\blue{\textcolor{blue}}
\def\printathome{\def\red{\null}\def\blue{\null}}

\date{\today}
\title[Decomposable polytopes with few vertices]{Decomposable polytopes with few vertices}

\author{Jie Wang}
\author{David Yost}

\address{Federation University, Mt. Helen, Vic. 3350, Australia}
\email{\texttt{wangjiebulingbuling@foxmail.com},\texttt{d.yost@federation.edu.au}}

\begin{abstract}
    A polytope $P$ is decomposable if it can be written as the Minkowski sum of two polytopes that are not homothetic to it. It is well known that decomposable $d$-dimensional polytopes with $2d$ or fewer vertices are the simplicial prisms. More recent work shows that, with two exceptions in low dimensions, decomposable $d$-dimensional polytopes with exactly $2d+1$ vertices are either pentasms or capped prisms. All these examples have a line segment for a summand. In this paper, we extend this conclusion to decomposable $d$-polytopes  with exactly $2d+2$ or  $2d+3$ vertices. In both cases, all such polytopes, with three exceptions,  have a line segment for a summand.  The exceptional cases are specified in detail.

\end{abstract}

\maketitle

\section{Introduction}

A polytope is the convex hull of a finite number of points in Euclidean space. A polytope $P$ is \emph{decomposable} if it can be written as the Minkowski sum of two polytopes, i.e. $P=Q+R=\{q+r: q \in Q, r \in R \}$, where $Q, R$ are not homothetic to $P$, i.e. not of the form $kP+a$ for $k>0$, $a\in \mathbb{R}^d$. Otherwise, it is called \emph{indecomposable}. 

The collection of all the faces of a polytope, ordered by inclusion, is a lattice, known as the \emph{face lattice}. Two polytopes $P$ and $Q$ are \emph{combinatorially equivalent} if their face lattices are isomorphic. Recall that a polytope is {\it combinatorially decomposable} if every polytope combinatorially equivalent to it is decomposable; a polytope is {\it combinatorially indecomposable} if every polytope combinatorially equivalent to it is indecomposable; a polytope is {\it conditionally decomposable} if one polytope combinatorially equivalent to it is decomposable and another one combinatorially equivalent to it is indecomposable.

For $d$-dimensional polytopes (in short, $d$-polytopes), the following are known:
\begin{itemize}
    \item All the $d$-polytopes with up to $2d-1$ vertices are indecomposable \cite[Lemma~5]{Yost1991};
    \item The only decomposable $d$-polytope with exactly $2d$ vertices is the simplicial $d$-prism \cite{Kallay1979}, \cite{Przeslawski-Yost-2016};
    \item Any decomposable $d$-polytopes with $2d+1$ vertices is either a pentasm, a 1-capped prism, $\Sigma(3)$, or $\Delta(2,2)$; all of them are combinatorially decomposable \cite[Theorem 5.9]{PUY2018Excess}. (Capped prisms, $\Sigma(3)$ and $\Delta(2,2)$ will be defined shortly.)  
\end{itemize}

In particular, there is no conditionally decomposable $d$-polytope with $2d+1$ or fewer vertices.

A natural extension of this topic is to characterise decomposable polytopes with a few more vertices. Our purpose is to give the classification of decomposable $d$-polytopes with $2d+2$, $2d+3$ vertices, respectively, and thereby to study the conditional decomposability of polytopes with few vertices. 

In dimension 3, there are exactly three conditionally decomposable polytopes with 8 vertices \cite{Yost2007}, which are $F282, F288, F287$ of Federico's list \cite{Federico1975}, where $Fn$ refers to the $n$th example in this catalogue. Those are the 3-dimensional conditionally decomposable polytopes with the minimum number of vertices. 
The first examples of conditionally decomposable polytopes in dimensions 4 and higher appeared in \cite{WJDYconditional}. For most of these conditional examples, the decomposable version has a line segment for a summand. This gave a hint to the importance of this property also for decomposability of polytopes with a small number of vertices. 

This work helps to classify the combinatorial types of decomposabe poyltopes with a small number of vertices, particularly in low dimensions. Note that we still do not know whether there is a conditionally decomposable 4-polytope with 11 vertices; some partial results are obtained in \cite{WJDYconditional}.

We need to know some sufficient conditions for the indecomposability of polytopes. We summarize the results which we need. Kallay \cite{Kallay1982} defined a {\it geometric graph} as any graph whose vertices lie in some Euclidean space and whose edges are some of the segments between them. This includes the \emph{edge graph}, i.e. the $1$-skeleton, of any polytope. He then generalized to such graphs the notion of decomposability of polytopes. We will say that a geometric graph $(V,E)$ is \emph{decomposable} if there is a {\it decomposing function} $\phi:V \to \mathbb{R}^d$ which is not the restriction of a homothety in $\mathbb{R}^d$. By a decomposing function we mean any function on the vertex set which respects the direction of edges, i.e. for any edge $e=[p,q]\in E$, there  is a scalar $\lambda_e$  such that $\phi(p)-\phi(q)=\lambda_e (p-q)$. Note that the scalar does not necessarily have to be positive. This presentation is not quite the same as Kallay's but it is equivalent. This generalisation is consistent, i.e. a polytope $P$ is decomposable if and only if its $1$-skeleton $G(V(P), E(P))$ is decomposable, where $V(P)$ stands for the set of vertices of $P$ and $E(P)$ stands for the set of edges of $P$ \cite[Corollary 5]{Kallay1982}. If $X$ is a set of vertices in graph, we denote by $G[X]$ the subgraph induced by $X$.

We say a facet $F$ of a polytope has \emph{Shephard's property} if every vertex in $F$ has a unique neighbor outside $F$. And we say a vertex $x$ is a \emph{simple extension} of $G$ if $G$ is an indecomposable graph, and $x$ is adjacent to at least two vertices in $G$.  

Before introducing our basic examples, we recall the notion of truncation. Let $H$ be a hyperplane intersecting the interior of $P$ and containing no vertex of $P$, and let $H^+$ and $H^-$ be  the two closed half-spaces bounded by $H$. Set $P':=H^+\cap P$, and suppose that $F$ is a face $P$ lying in $H^-$. Then the polytope $P'$ is said to be obtained by {\it truncating} the face $F$ by $H$. 

Denote by $\Delta(m,n)$ the Cartesian product of a simplex in $\mathbb{R}^m$ and a simplex in $\mathbb{R}^n$, or any polytope combinatorially equivalent to it. It is a simple $(m+n)$-dimensional polytope with $(m+1)(n+1)$ vertices and $m+n+2$ facets. In particular, $\Delta(d-1,1)$ is the Cartesian product of a $(d-1)$-simplex and a line segment; we will refer to it henceforth simply as a {\it prism}, or a $d$-prism if the dimension needs to be specified.

Denote by $M(k,d-k)$ a $(d-k)$-fold pyramid over a $k$-prism. 
Denote by $C(d)$ a $d$-polytope obtained by truncating a simple edge of $M(2,d-2)$. 
Denote by $\Sigma(d)$ a certain polytope which can be expressed as the Minkowski sum of two simplices, combinatorially equivalent to $\text{conv}\{ 0, e_1, e_1+e_k, e_2, e_2+e_k, e_1+e_2, e_1+e_2+2e_k: 3\le k\le d  \}$, where $e_i$ is the standard basis of $\mathbb{R}^d$. 

Denote by \emph{$d$-pentasm} the Minkowski sum of a simplex and a line segment which is parallel to one triangular face of the simplex but not parallel to any edge of the simplex, or any polytope combinatorially equivalent to it. Alternatively, it can be seen as the polytope obtained by truncating a simple vertex of $M(2,d-2)$. 
The \emph{$d$-dimensional 1-capped prism} is obtained by stacking a simplex (a pyramid) over one of the simplex facets of the prism, where the apex may or may not lie in some 2-faces of the prism. Here, stacking means gluing a simplex onto one of simplex facets of the polytope. There are $d-2$ different combinatorial types of a capped prism; for more details, see \cite[p. 2015]{PUY2018Excess}.

The followings are the necessary material that we will need.

\begin{theorem}\label{lit}
	Let $P$ be a $d$-polytope.
	\begin{enumerate}
		\item If $P$ is decomposable, then it has at least $d$ decomposable facets \cite[Theorem 4.4]{PUY2018Excess}. 
		\item If there is an indecomposable subgraph of $P$ that touches every facet of $P$, then $P$ is indecomposable \cite[Theorem 9]{Przeslawski-Yost-2008}.
		\item If $P$ has a facet with Shephard's property, then $P$ is combinatorially decomposable \cite[result 15]{Shephard1963}.
		\item  A pyramid is indecomposable.
		\item If $x$ is a simple extension of an indecomposable graph $G$, then $G \cup \{x\}$ is an indecomposable graph \cite[Proposition 1]{Przeslawski-Yost-2016}.
		\item If $G_1$, $G_2$ are indecomposable subgraphs,and they are connected by skew lines, then $G_1 \cup G_2$ is indecomposable \cite[Theorem 10]{Kallay1982}. 
		\item If $P$ is a decomposable polytope, $F$ is a facet of $P$, and the subgraph of $G(P)$ disjoint from $F$ is indecomposable, then $F$ has Shephard's property. 
		\item The subgraph induced by all the vertices outside a facet $F$ (of a polytope) is connected \cite[Theorem 15.5]{Brondstedbook}. 
	\end{enumerate}
\end{theorem}

We will freely use the fact that {\bf projective transformations do not affect decomposability}. That is, if $P$ is any polytope and $T$ is a projective transformation, then $P$ and $T(P)$ are either both decomposable or both indecomposable \cite{Kallay1984, Smilansky1987}. 

Throughout the paper, when we say a polytope $P$ has a line segment for a summand, we mean that there is a projective transformation $T$ such that $T(P)$ has a line segment for a summand. Note that $T(P)$ is combinatorially equivalent to $P$. We say two polytopes $P,Q$ are \emph{projectively equivalent} if there is a projective transformation $T$ such that $T(P)=Q$. Similarly, when we say \emph{projectively parallel} edges $\{ e_1, \dots, e_i  \}$, we mean that there is a projective transformation $T$ such that the edges $\{ T(e_1), \dots, T(e_i) \} \subset T(P) $ are all parallel, where the polytope $T(P)$ is combinatorially equivalent to $P$.

\begin{remark}\label{cross-section}
Suppose that a polytope $P$ has a line segment for a summand, i.e. $P=Q+I$, where $I=[a, b]$ denotes a line segment. Since any vertex of a Minkowski sum of polytopes decomposes uniquely as the sum of vertices from the respective summands, see \cite[p. 317, {\bf1}]{Grunbaum2003}, we then have a natural partition of the vertices of $P$ according to which vertex of $I$ they come from. This defines a decomposing function $f:V(P)\to\{a,b\}$, with $f^{-1}(a)$ and $f^{-1}(b)$ being a partition of $V(P)$. For convenience,  the edges that are parallel to $[0,a]$ will be called ``side" edges of $P$ and the faces that contain some of the ``side" edges will be called ``side" faces of $P$. It follows that there is a polytope $P'$ combinatorially equivalent to $P$ such that there is a hyperplane (i.e. a cross section) that passes through all the ``side" edges without passing through any vertex, which separates $V(P)$ into two parts, namely $f^{-1}(a)$ and $f^{-1}(b)$.
\end{remark}

\begin{figure}[H]
	     \centering
	     \includegraphics[width=0.25\textwidth]{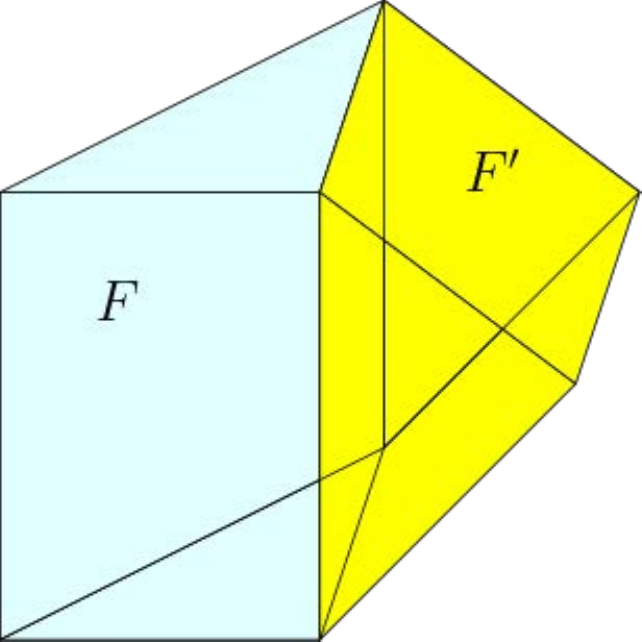}
	     \caption{``rotated" prisms}
	     \label{fig:rotated_illustration}
     \end{figure}

One of the main observation in this paper is that in a decomposable $d$-polytope with $2d+2$ or $2d+3$ vertices, when the dimension is high (say $\ge 6$), a decomposable facet $F$ must have a line segment for a summand and the other facet $F'$ that intersects with $F$ at a ``side" ridge must be decomposable, and then $F'$ also has a line segment for a summand. Therefore, there is at least an edge $e\in F'$ outside $F$ that is projectively parallel to the ``side" edges of $F$ and we shall say that the edge $e$ is \emph{associated} with $F'$. In lower dimensions, things get tricky, especially in dimension 4, where the ``directions" of those sets of projectively parallel edges do not have to be coincide. In other words, we can have two 3-prism facets $F, F'$ which intersect in a quadrilateral, and if we look at the two vertices that is outside the quadrilateral and the edge incident to the two vertices, we can see that these two edges can have different ``directions". Then we would say that the other facet $F'$ is \emph{``rotated"}; see Figure~\ref{fig:rotated_illustration}.

\section{Decomposable $d$-polytopes with $2d+2$ vertices, $d\ge 4$}



We first discuss decomposable 4-polytopes with 10 vertices. 

\begin{lemma}
	Let $P$ be a decomposable $4$-polytope with $10$ vertices. Then there is a decomposable facet $F$ with maximum number of vertices, which has $6$, $7$ or $8$ vertices.
	\begin{proof}
		Since a decomposable $d$-polytope has at least $d$ decomposable facets (by \cref{lit}(1)), we may let $F$ be a decomposable facet with the maximum number of vertices. As $F$ is a 3-dimensional decomposable polytope, it has at least 6 vertices. As $P$ is decomposable, it cannot be a pyramid, so $F$ has at most 8 vertices.
	\end{proof}
\end{lemma}

Recall that if $f$ is a decomposing function of $P$, and each scalar $\alpha$ is in the interval $[0,1]$, then two summands of $P$ are the convex hull of $f(V(P))$ and the convex hull of $(id_f)(V(P))$, where $id$ stands for the identity map.  

The following lemma is embedded in the proof of Theorem 5 of \cite{WJDYconditional}.

\begin{lemma}\label{2verticesoutside}
	Let $F$ be a facet of a decomposable polytope $P$ such that there are exactly two vertices outside $F$. Then $F$ has a line segment for a summand and $P$ has a line segment for a summand.
\end{lemma}

The same reasoning in \cref{2verticesoutside} gives the following.

\begin{lemma}\label{decomposable:triangleoutside}
	If $F$ is a decomposable facet of a decomposable polytope and there are three vertices pairwise adjacent that are outside of $F$, then $F$ and $P$ both have a triangle for a summand.
\end{lemma}

\begin{corollary}\label{dim4:linesegment}
	Let $P$ be a decomposable $4$-polytope with $10$ vertices, and let $F$ be a decomposable facet with $8$ vertices. Then $P$ has a line segment for a summand.
\end{corollary}

A suitable partition of the vertices may give us a clue to its summands.

\begin{lemma}\label{decompfunc-linesegment}
	Let $V_1, V_2$ be the partition of vertices of $P$, such that each vertex in $V_1$ (or $V_2$, resp.) has at most one neighbor in $V_2$ (or $V_1$, resp.). If the edges in between $V_1$ and $V_2$ are projectively parallel, then $P$ has a line segment for a summand.
	\begin{proof}
		By taking a suitable projective transformation, we may assume that the edges between $V_1, V_2$ in $P$ are all parallel. We can define a decomposing function $f$ on $P$, where for $v\in V_2$, $f(v)=v$, and for $w\in V_1$, $f(w)=w-a$, where $a$ is the direction of the parallel edges between $V_1$ and $V_2$, and $\|a \| $ is sufficiently close to 0. The result now follows.
	\end{proof}
\end{lemma}

\begin{lemma}\label{lem:cappedpentasm_notrianglesummand}
	A capped prism or pentasm does not have a triangle for a summand.	
	\begin{proof}
		Let $V_1$, $V_2$ be the set of vertices such that the edges between them are the projectively parallel edges. Note that $G[V_1], G[V_2]$ are indecomposable. Let $f$ be a decomposing function on a capped prism or pentasm. Then $f(V_i)$ is either a single point or a homothety of $V_i, i=1,2$. So, the image of $f$ on all the vertices cannot be a triangle.		
	\end{proof}
\end{lemma}

\begin{lemma}\label{lem:Delta22_nolinesegment}
	The polytope $\Delta(2,2)$ does not have a line segment for a summand.	
	\begin{proof}
		Let $V_1, V_2, V_3$ be the set of vertices of the triangular faces of $\Delta(2,2)$ that form a partition of vertices of $\Delta(2,2)$. Let $f$ be a decomposing function on $\Delta(2,2)$. Then $f(V_i)$ is either a single point or a homothety of $V_i, i=1,2,3$. Suppose $\Delta(2,2)$ has a line segment $[a,b]$. Then $f(V_i)$ is  either $a$ or $b$, where $i=1,2,3$. By the pigeonhole principle, (say) $V_1$ and $V_2$ are mapped to the same point. This is not possible, as edges between $V_1,V_3$ and edges between $V_2, V_3$ cannot all be parallel to $[a,b]$.
	\end{proof}
\end{lemma}

\begin{lemma}\label{dim4:sigma3}
	Let $P$ be a decomposable $4$-polytope with $10$ vertices. If one of the decomposable facets with the maximum number of vertices is $\Sigma(3)$, then $P$ is either $\Sigma(4)$ or $C(4)$.
	\begin{proof}
		Let $F$ be a facet which is $\Sigma(3)$. If the three vertices outside $F$ are pairwise adjacent, then, by \cref{lit}(7), $F$ has Shephard's property. In this case, $P$ has 21 edges. Therefore $P$ is either $C(4)$ or $\Sigma(4)$ \cite[Theorem 6.1]{PUY2018Excess}. It follows that $P$ is $C(4)$ (as in $\Sigma(4)$, the three vertices outside any $\Sigma(3)$ facet are not pairwise adjacent).
		
		Suppose now the three vertices outside $F$ are not pairwise adjacent, label them as $x_1,x_2$ and $y$; see Figure~\ref{fig:dim4_sigma3}. Consider the facets $T_1$, $T_2$ intersecting $F$ at a triangular ridge. Then $T_1$ and $T_2$ cannot be a prism, as $x_1$ is not adjacent to $x_2$. So if $T_1$ or $T_2$ has 6 vertices, it must be indecomposable, then $P$ would be indecomposable, by \cref{lit}(2). Let $v_4 \in F$ be the common vertex of $T_1$ and $T_2$. If $T_1$ or $T_2$ has 5 vertices, then either $T_1$ and $T_2$ share two common vertices, or their graphs are connected by a triangle, or $x_1$-$y$-$x_2$-$v_4$-$x_1$ forms a non-planar 4-cycle, in each case $P$ is indecomposable. 
		
		Thus $T_1$ and $T_2$ can only be simplices. Note that the facets $T_1$ and $T_2$ have different apexes and $y$ is not an apex of them. Suppose the contrary that $T_1, T_2$ have the same apex, then $G[T_1 \cup T_2]$ is indecomposable. It follows that the other two vertices of $F$ are simple extensions of $G[T_1 \cup T_2]$, which leaves out 2 vertices, making $P$ indecomposable. Assume now without loss of generality that $y$ is an apex of $T_1$ and $x_2$ is an apex of $T_2$. Then $G[v_4 \cup y \cup x_2]$ is a triangle (indecomposable), which shares an edge with $T_1$ and shares an edge with $T_2$. It follows that $G[T_1 \cup T_2 \cup x_2 \cup y]$ is indecomposable. As the other two vertices of $F$ are simple extensions of $G[T_1 \cup T_2 \cup x_2 \cup y]$, the polytope $P$ would be indecomposable. 		
		So one of the simplices has apex $x_1$, the other one has apex $x_2$. Without loss of generality, we assume that $x_1$ is the apex of $T_1$, and $x_2$ is the apex of $T_2$.
		\begin{figure}[H]
			\centering
			\includegraphics[width=0.45\textwidth]{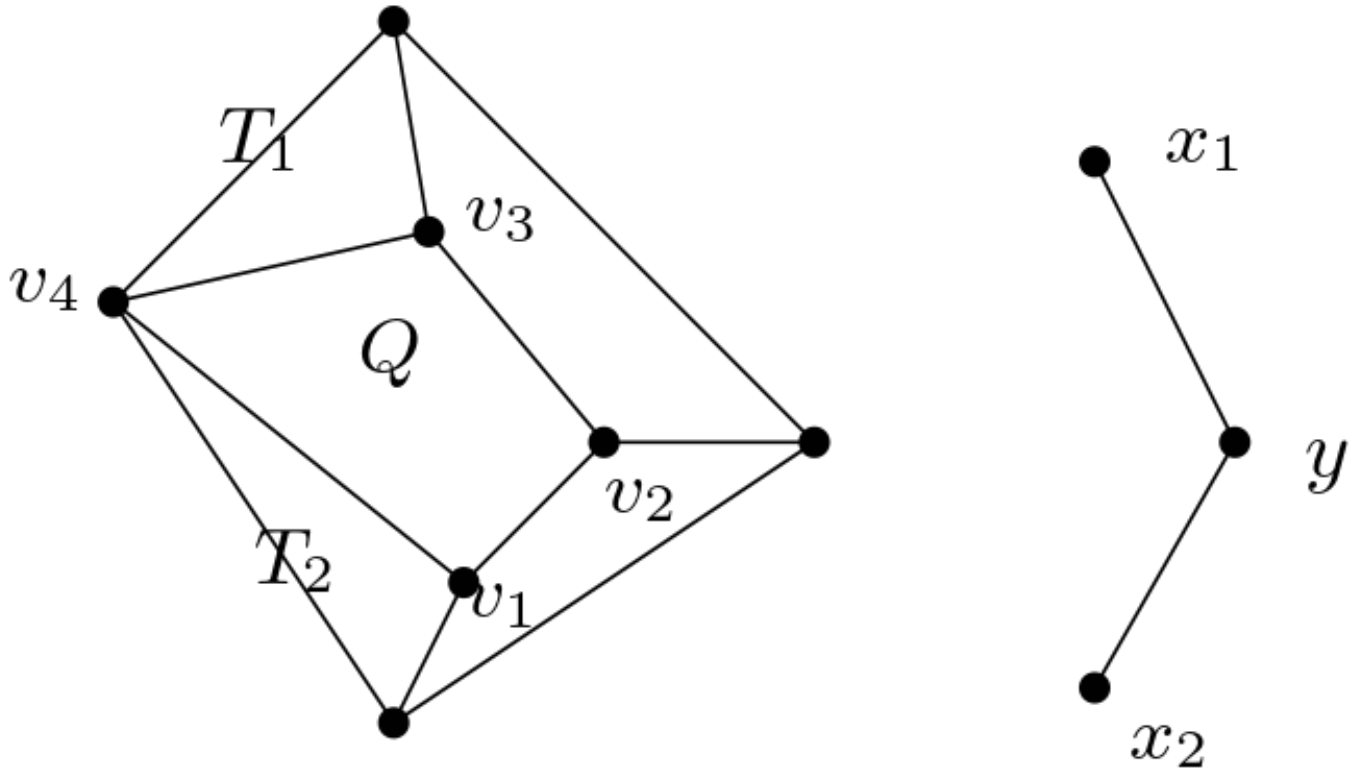}
			\caption{three vertices outside $F$ are not pairwise adjacent}
			\label{fig:dim4_sigma3}
		\end{figure}
		The facets that intersect $F$ at a quadrilateral ridge cannot be indecomposable. Otherwise, by \cref{lit}(2), $P$ would be indecomposable. Let $Q$ be a facet that intersects $F$ at a quadrilateral ridge which intersects both triangular faces of $F$; see Figure~\ref{fig:dim4_sigma3}. Then $Q$ can be $\Sigma(3)$, a pentasm, a capped prism or a prism. 
		
		If $Q$ is a prism, then $y \in Q$. By symmetry, there are only two cases to consider: $y$ is adjacent to $v_1, v_2$ and $y$ is adjacent to $v_1, v_4$. In the former case, $G[y \cup v_1 \cup v_2 \cup T_2]$ is indecomposable (simply by a lot of triangles). Then $x_1$ serves as a simple extension of $G[y \cup v_1 \cup v_2 \cup T_2]$. It follows that $G[y \cup v_1 \cup T_1 \cup T_2]$ is indecomposable, implying that the polytope is indecomposable (by \cref{lit}(2)). In the latter case, $G[y \cup T_2]$ is indecomposable. And again, $x_1$ serves as a simple extension of $G[y \cup T_2]$ and so, $G[y \cup T_1 \cup T_2]$ is indecomposable, implying that the polytope is indecomposable, a contradiction. 
		
		Since $x_1$ and $x_2$ are both adjacent to $v_4$, $Q$ cannot be a pentasm or a capped prism, therefore, $Q$ is $\Sigma(3)$. This completes the entire edge graph of $P$, and $P$ must be $\Sigma(4)$ (by \cite[Theorem 4.9]{PUY2020} that a polytope with excess degree at most $d-1$ is determined by its edge graph).
	\end{proof}
\end{lemma}

\begin{lemma}\label{dim4:nonesigma3}
	Let $P$ be a decomposable $4$-polytope with $10$ vertices. If one of the decomposable facets $F$ with the maximum number of vertices is a pentasm or a capped prism and none of the facets of $P$ is $\Sigma(3)$, then $P$ has a line segment for a summand.
	\begin{proof}
		By \cref{lem:cappedpentasm_notrianglesummand}, pentasms or capped prisms do not have a triangle for a summand, so by \cref{decomposable:triangleoutside}, the three vertices outside $F$ are not pairwise adjacent. Let $F'$ be the other facet that intersects $F$ at a quadrilateral/pentagon. We claim that $F'$ is decomposable. Suppose that $F'$ is indecomposable. Then the subgraph of $F \cup F'$ is indecomposable (by simple extensions), which leaves out at most 2 vertices of $P$. In this case, by \cref{lit}(2), $P$ would be indecomposable. Let $U, V$ be the partition of vertices of $F$ such that the edges in between $U, V$ are projectively parallel. Note that the subgraphs induced by $U$, $V$ are both indecomposable subgraphs.

		Denote the three vertices outside $F$ as $x_1, x_2, x_3$, and assume, without loss of generality, that $x_1$ is not adjacent to $x_3$. Since deg($x_1$)$\ge 4$,  we know that $x_1$ is adjacent to at least 3 vertices in $F$. 
		
		We claim that $x_1$ is only adjacent to either vertices of $U$ or vertices of $V$. Suppose the contrary that $x_1$ is adjacent to at least 2 vertices in $U$, and at least 1 vertex in $V$. Since there is always a path of length 2 from a vertex in $U$ to a vertex in $V$, this guarantees a 4-cycle containing $x_1$, one neighbor $u\in U$ of $x_1$, one neighbor $v\in V$ of $x_1$, and the other vertex that connects $u, v$. It is not hard to see that this is a non-planar 4-cycle. So the subgraph $G[x_1 \cup U]$ is indecomposable. Note that the neighbor $v\in V$ of $x_1$ is a simple extension of $x_1 \cup U$, so the subgraph $G[x_1 \cup v \cup U]$ is indecomposable, then the subgraph $G[x_1 \cup V \cup U]$ is indecomposable (again by simple extensions). So, by \cref{lit}(2), $P$ would be indecomposable, a contradiction. The same is true for $x_3$. 
		
		\begin{figure}[H]
				\centering
				\includegraphics[width=\textwidth]{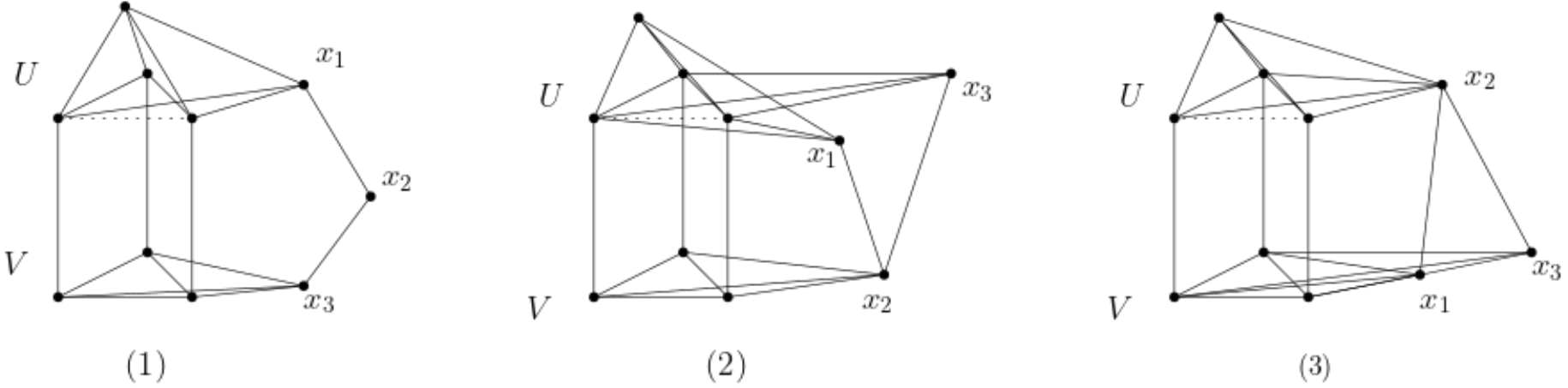}
		    	\caption{$F$ has 7 vertices, and the three vertices outside $F$ are not pairwise adjacent}
		    	\label{fig:Fhas7}
		\end{figure}
		If $x_1, x_3$ are both only adjacent to vertices in $U$, then $x_2$ is adjacent to all the vertices in $V$, see Figure \ref{fig:Fhas7}(2). The edges $[x_1, x_2]$ and $[x_2, x_3]$ cannot both be projectively parallel to the ``side" edges of $F$. Then we have two indecomposable subgraphs connected by skew lines. By \cref{lit}(6), $P$ is indecomposable, a contradiction. The same for the case that $x_1, x_3$ are both adjacent  only to vertices in $V$. 
		
		Suppose now that $x_1$ is only adjacent to $U$, and $x_3$ is only adjacent to $V$; see Figure~ \ref{fig:Fhas7}(3). We will show that $x_2$ is adjacent either only to vertices in $U$ or only to vertices in $V$. Suppose $x_2$ is adjacent to a vertex in $U$. Then $x_2$ is a simple extension of the subgraph $G[x_1 \cup U]$, so now we have two indecomposable subgraphs $x_1 \cup x_2 \cup U$ and $x_3 \cup V$. Additionally, any edge in between these two indecomposable subgraphs has to be projectively parallel to the ``side" edge of $F$. It follows that $x_2$ is only adjacent to vertices in $U$, and the lines containing $[x_2, x_3]$ are projectively parallel to the lines containing the ``side" edges of $F$. By \cref{decompfunc-linesegment}, $P$ has a line segment for a summand.	
	\end{proof}
\end{lemma}

\begin{lemma}\label{dim4:allprisms}
	Let $P$ be a decomposable $4$-polytope with $10$ vertices. If all of its decomposable facets are triangular prisms, then $P$ has a line segment for a summand.	
	\begin{proof}
		Let $F$ be a prism facet, and let $F'$ be the other facet that intersects with $F$ at a quadrilateral ridge (of $P$). We claim that $F'$ is decomposable. Suppose that $F'$ is indecomposable. Then the subgraph of $F \cup F'$ is indecomposable, which leaves out at most 3 vertices of $P$. It follows that, by \cref{lit}(2), $P$ would be indecomposable. So every other facet $F'$ that intersects $F$ at a quadrilateral is a prism. Without loss of generality, we may denote the sets of vertices of the triangles of $F$ by $U, V$ and may assume that edges in between $U, V$ are parallel. The following are the possible subgraphs induced by the 4 vertices outside $F$; see Figure~\ref{fig:prism}.		
		
		\begin{figure}[H]
			\centering
			\includegraphics[width=0.95\textwidth]{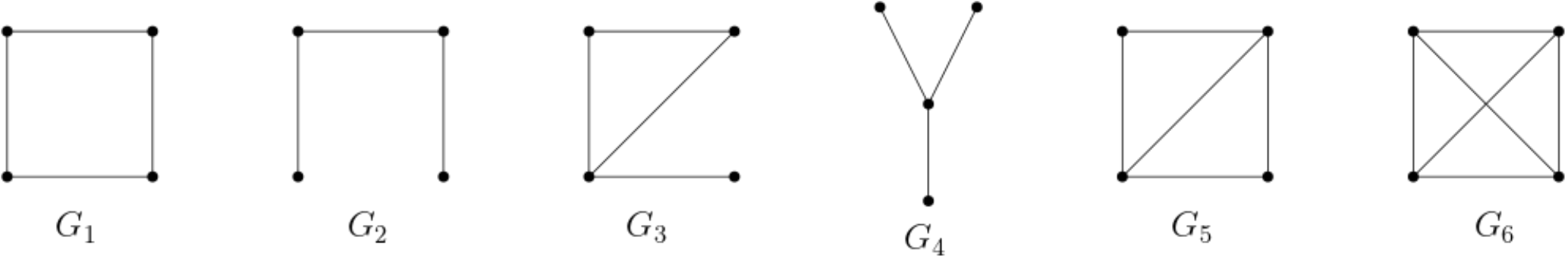}
				\caption{The subgraph induced by the 4 vertices outside $F$}
				\label{fig:prism}
		\end{figure}
	
	    We deal with each graph in Figure~\ref{fig:prism} at a time.
	
	    \textcolor{blue}{$G_1$}. 
	    
	    Suppose $F'$ is a facet that intersects $F$ at a quadrilateral and the parallel edges of $F'$ are in the different direction, namely, $F'$ is ``rotated". Assume without loss of generality that $x_1, x_2 \in F'$; see Figure~\ref{fig:rotated}.
	    \begin{figure}[H]
	   	   \centering
		   \includegraphics[width=0.35\textwidth]{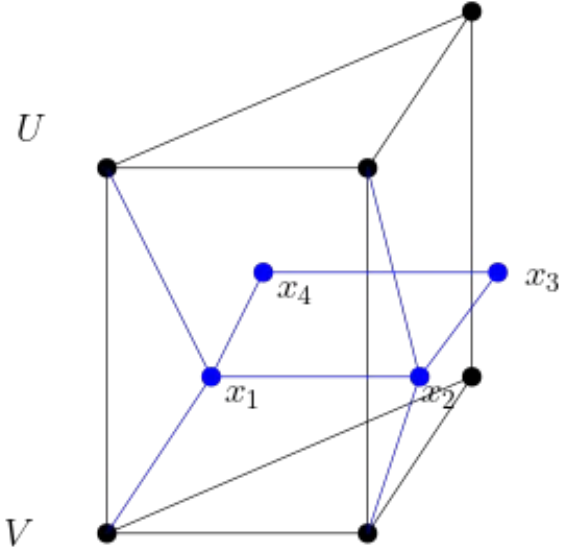}
		   	\caption{$F'$ is ``rotated"}
		    \label{fig:rotated}
       	\end{figure}
		In this scenario, $x_1$ cannot be adjacent to any other vertices except $x_4$. Moreover, $x_2$ cannot be adjacent to any other vertices except $x_3$. Suppose otherwise that $x_1$ is adjacent to another vertex in $U$. Then the subgraph $G[x \cup U]$ is indecomposable. Since now $x_2$ is a simple extension of $G[x \cup U]$, we have $G[x_1 \cup x_2 \cup U]$ indecomposable. Then $G[x_1 \cup x_2 \cup U \cup V]$ is indecomposable (because of the triangular faces in $F'$ and simple extension of the third vertex in $V$). By \cref{lit}(2), $P$ would be indecomposable, a contradiction. The same applies to the case that $x_1$ is adjacent to another vertex in $V$. So, $x_1, x_2$ are not adjacent to any other vertex except $x_4, x_3$, respectively. It follows that $x_1, x_2$ are simple vertices.  
		
		Consider another facet $F''$ that intersects $F$ at a quadrilateral face. Suppose $F''$ is not ``rotated". Note that $x_2 \notin F''$, as $x_2$ is adjacent to vertices both in $U$ and $V$. Since $x_1$ is not adjacent to any other vertices in $U$ and $V$, we know that $x_1 \notin F''$ either. So, we must have $x_3, x_4 \in F''$ and the edge $[x_3, x_4]$ is projectively parallel to the ``side" edges of $F$. We may assume without loss of generality that $x_3$ is adjacent to two vertices in $V$ and $x_4$ is adjacent to two vertices of $U$. By looking at the triangles $\triangle_{v_1, v_2, v_3}, \triangle_{v_2, v_3, x_3}, \triangle_{x_2, x_3, v_3}$, we can see that the subgraph $G[u_3 \cup x_2 \cup x_3 \cup V]$ is indecomposable; see Figure~\ref{fig:rotated2}. It follows that $x_1, x_4$ are simple extensions of $G[u_3 \cup x_2 \cup x_3 \cup V]$, making $P$ indecomposable, a contradiction. 		
			\begin{figure}[H]
			  \centering
			  \includegraphics[width=0.35\textwidth]{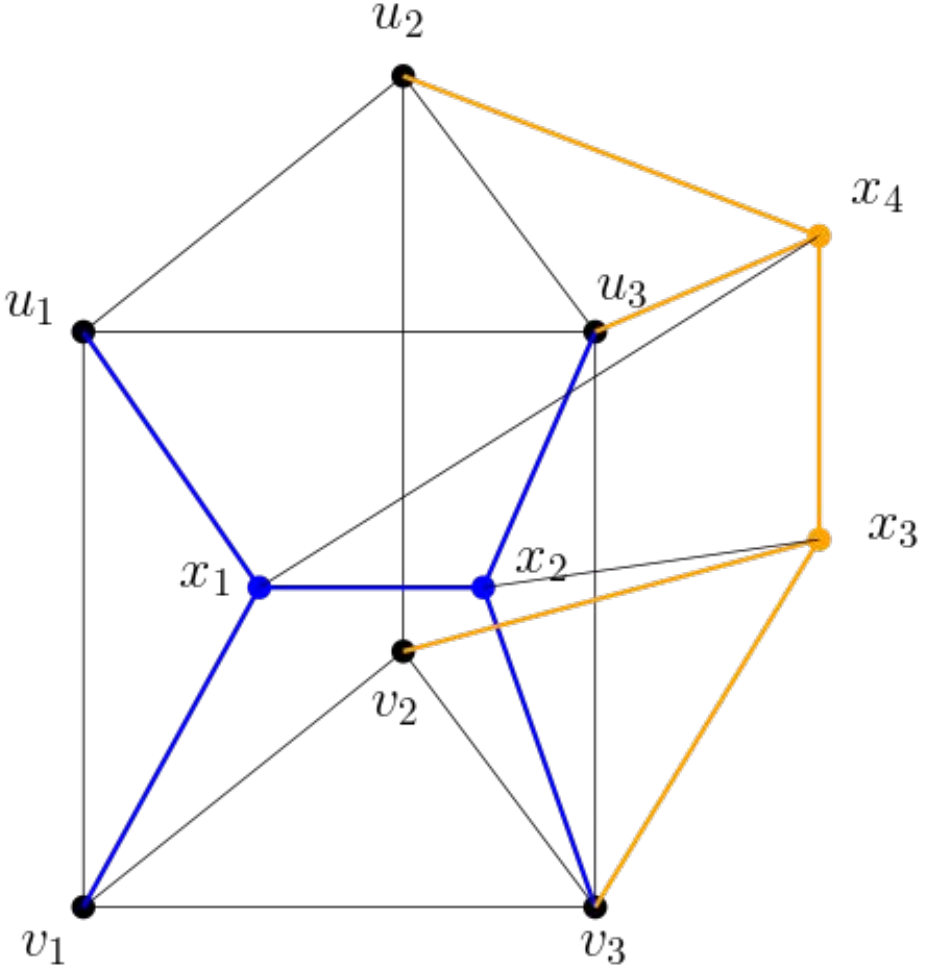}
			 	\caption{$F'$ is ``rotated" but $F''$ is not}
			    \label{fig:rotated2}
	     	\end{figure}
    	Now for every other facet $F''$, that intersects $F$ at a quadrilateral, $F''$ is ``rotated". By the previous discussion, we know that $x_1, x_2, x_3, x_4$ are all simple. So $P$ has $9+8+4=21$ edges, by \cite[Theorem 6.1]{PUY2018Excess}. There are only four 4-polytopes with 10 vertices, with 21 or fewer edges; those are $A(4), B(4), C(4)$ and $\Sigma(4)$. However, they all have a decomposable facet with at least 9 vertices, so, those are not possible.
    	
    	So, the only possible case is that every other facet $F'$ which intersects $F$ at a quadrilateral face is not ``rotated". We may assume without loss of generality that $[x_1, x_2]$ is the edge associated to $u_1, v_1, u_3, v_3$ which together with $[x_1, x_2]$ forms a facet; see Figure~\ref{fig:allprism_nice}. Note that $x_1$ cannot be adjacent to any vertex in $V$, and $x_2$ cannot be adjacent to any vertex in $U$, otherwise we have skew lines connecting two indecomposable subgraphs that leave out only two vertices outside, making $P$ indecomposable. Consider the neighbors of $x_4$: if $x_4$ is adjacent to a vertex in $U$, then the subgraph $G[x_1 \cup x_4 \cup U]$ is indecomposable, so $x_4$ cannot be adjacent to any vertex in $V$, otherwise we have skew lines connecting two indecomposable subgraphs that leave out only two vertices outside; if $x_4$ is adjacent to two vertices (or more) in $V$, then $x_2 \cup x_4 \cup V$ is indecomposable, and the line containing $[x_1, x_4]$ is skew to one of the lines containing $[u_i, v_i], i=1,2,3$, then $P$ would be indecomposable. So, $x_4$ is only adjacent to vertices of $U$. Similarly, $x_3$ is only adjacent to vertices of $V$. Now we have a partition of vertices of the polytope, $T_1=\{x_1, x_4, U\}$, $T_2=\{x_2, x_3, V\}$. The subgraphs induced by $T_1, T_2$ are indecomposable, so, the edges in between them are projectively parallel. By \cref{decompfunc-linesegment}, $P$ has a line segment for a summand.
       	\begin{figure}[H]
        	\centering
        	\includegraphics[width=0.45\textwidth]{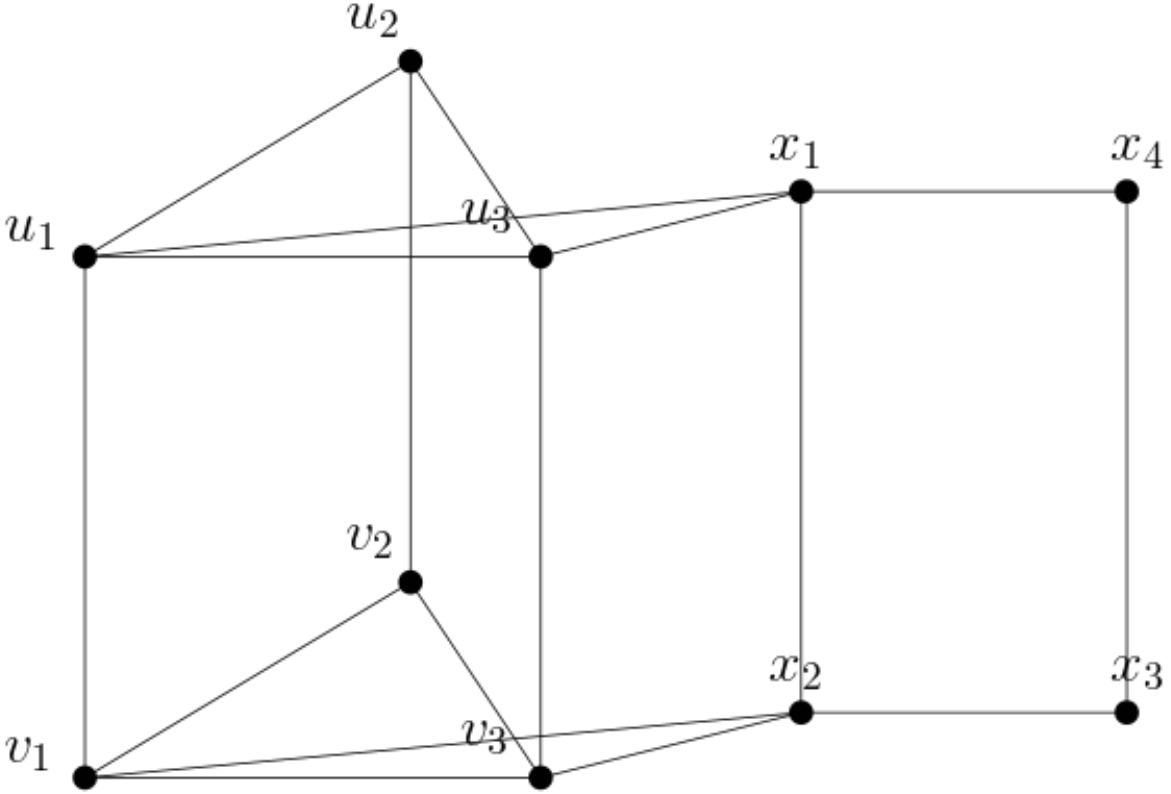}
        	\caption{all the other facet $F'$ are not ``rotated"}
        	\label{fig:allprism_nice}
        \end{figure}
    
    	\textcolor{blue}{$G_2$}. 
    	
    	Label the vertices of $G_2$ as $x_1, x_2, x_3, x_4$ such that $x_1, x_4$ are not adjacent. Since $P$ is a 4-polytope, $x_1, x_4$ have at least 3 neighbors in $F$. It is not hard to see that their neighbors are either only in $U$ or only in $V$, since otherwise, there is a skew line connecting two indecomposable subgraphs that leave out at most 3 vertices outside, implying that $P$ would be indecomposable (by \cref{lit}(2)). 
    	
		Consider the other facet $F'$ that intersects $F$ at a quadrilateral face. If $F'$ were ``rotated", then $x_2, x_3\in F'$, since $x_1, x_4$ cannot be adjacent to vertices both in $U$ and $V$. But then we have a collection of triangles; see Figure~\ref{fig:G2_1}. And by some simple extensions, one can see that this case cannot occur. So, every other facet that intersects $F$ at a quadrilateral face has ``side" edges projectively parallel to the ``side" edges of $F$.
		\begin{figure}[H]
			\centering
			\includegraphics[width=0.45\textwidth]{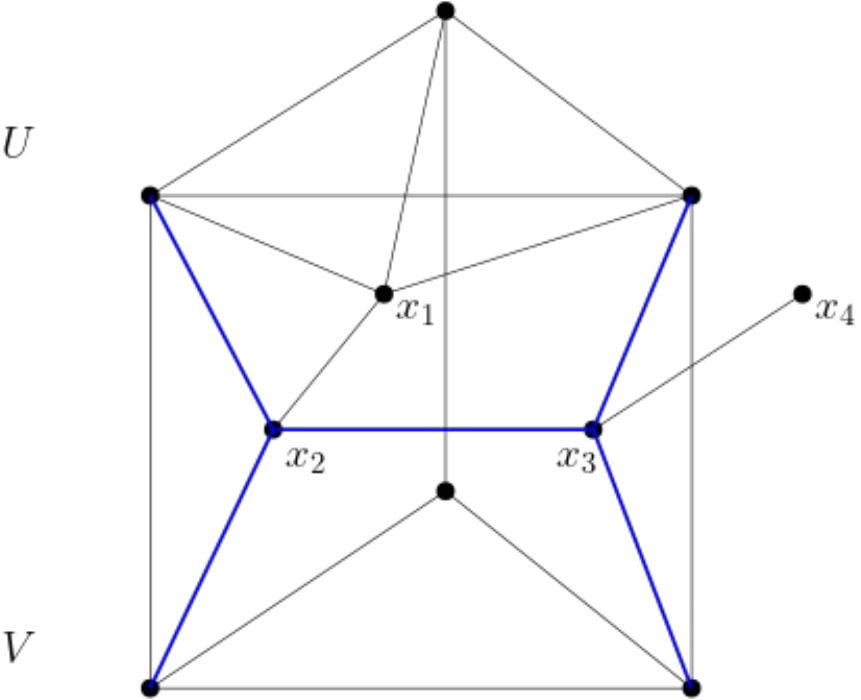}
			\caption{$G_2$, ``rotated" $F'$}
			\label{fig:G2_1}
		\end{figure}
	
	    Now it suffices to show that $x_2, x_3$ are not adjacent to both vertices in $U$ and vertices in $V$. We may assume that $x_1$ is adjacent to all the vertices in $U$. If $x_2$ is adjacent to $u\in U$ and some vertex in $v\in V$, the line containing $[x_2, v]$ is skew to one of the lines containing the projectively parallel ``side" edges of $F$. However, since the subgraph $G[x_1, x_2, U]$ is indecomposable, the subgraph induced by $V$ is indecomposable. This gives an indecomposable subgraph that leaves out 2 vertices outside, making $P$ indecomposable. The same is true for $x_3$. By the adjacency of $x_i$ and $U,V$, we can give the partition of the vertices, where the subgraphs of each part is indecomposable, so the edges connecting them must be projectively parallel. By \cref{decompfunc-linesegment}, $P$ has a line segment for a summand.  
	    
	    \textcolor{blue}{$G_3$}. 
	    
	    Let $x_1$ and $x_2$ be the vertices of $G_3$ such that $x_1$ has degree 1 and $x_2$ has degree 3. Again, we know that $x_1$ is adjacent either only to vertices in $U$ or only to vertices in $V$. We may assume without loss of generality that $x_1$ is only adjacent to $U$. 
	     
	    Suppose the other facet $F'$ that intersects $F$ at a quadrilateral face is a ``rotated" prism. We first show that $x_3, x_4 \notin F'$. Suppose the contrary that $x_3, x_4 \in F'$.  
	    
	    Case 1: If $x_2$ is adjacent to a vertex in $U$, then $x_1 \cup x_2 \cup U$ is indecomposable, and $x_3, x_4$ are simple extensions of this subgraph, so we will have an indecomposable subgraph that leaves only 3 vertices outside, making $P$ indecomposable; see Figure~\ref{fig:G3_1}.
	    \begin{figure}[H]
	    	\centering
	    	\includegraphics[width=0.35\textwidth]{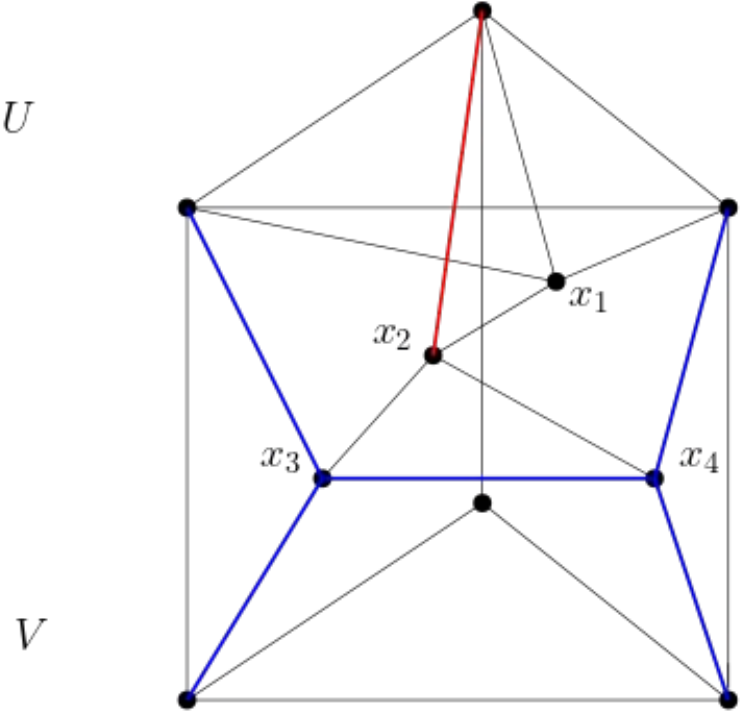}
	    	\caption{$G_3$, ``rotated" $F'$, case 1}
	    	\label{fig:G3_1}
	    \end{figure}
    
        Case 2: $x_2$ is adjacent to some vertex in $V$. If $x_2$ and $x_3$ (or $x_2$ and $x_4$, respectively) are adjacent to the same vertex in $V$, then we have an indecomposable subgraph $G[x_2 \cup x_3 \cup u_1 \cup v_1]$. This subgraph is made of triangles, which implies that $x_1, x_4$ are simple extensions of it. Thus, this indecomposable subgraph can be extended to every vertex, making $P$ indecomposable; see Figure~\ref{fig:G3_2_1}. The remaining case is that $x_2$ is only adjacent to $v_2$; see Figure~\ref{fig:G3_2_2}. Note that $x_3, x_4$ are simple (the argument is the same as for $G_1$). Consider the other facet $F''$ that intersects $F$ at $u_1,v_1, u_2, v_2$. Since $v_2$ is only adjacent to $x_2$ (as $x_1$ is not adjacent to any vertex in $V$, and $x_3, x_4$ are simple), $x_2 \in F''$. So either $[x_1, x_2]$ is the edge associated (the unique edge that is outside the facet) with $F''$, or $[x_2, x_3]$ is the edge associated with $F''$. But this is not possible. In the former case, $x_2$ is not adjacent to $v_1$, and in the latter case, $x_2$ is not adjacent to $u_2$, so, by the structure of a prism, this is not possible. 
        
        \begin{figure}[H]
        	\centering
        	\begin{minipage}{.5\textwidth}
        		\centering
        		\includegraphics[width=.6\linewidth]{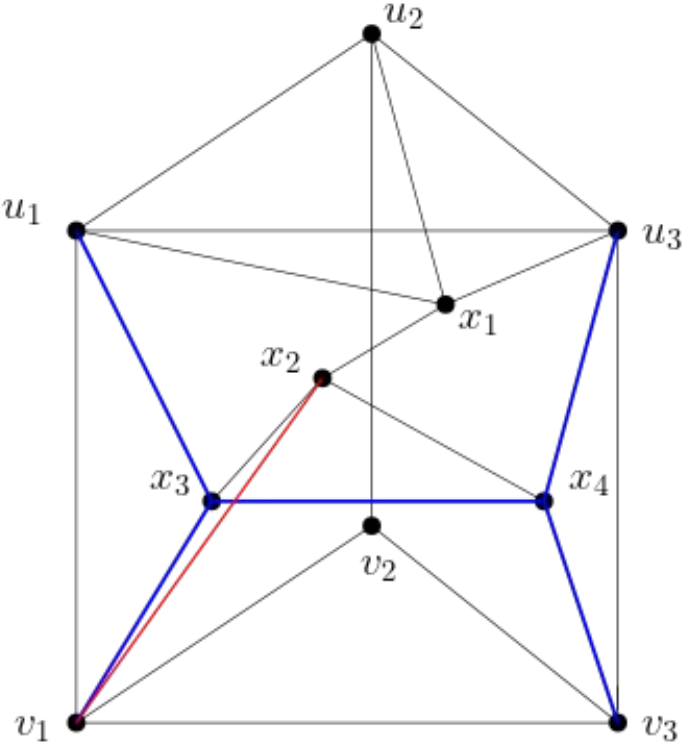}
        		\captionof{figure}{$G_3$, ``rotated" $F'$, case 1 }
        		\label{fig:G3_2_1}
        	\end{minipage}%
        	\begin{minipage}{.5\textwidth}
        		\centering
        		\includegraphics[width=.6\linewidth]{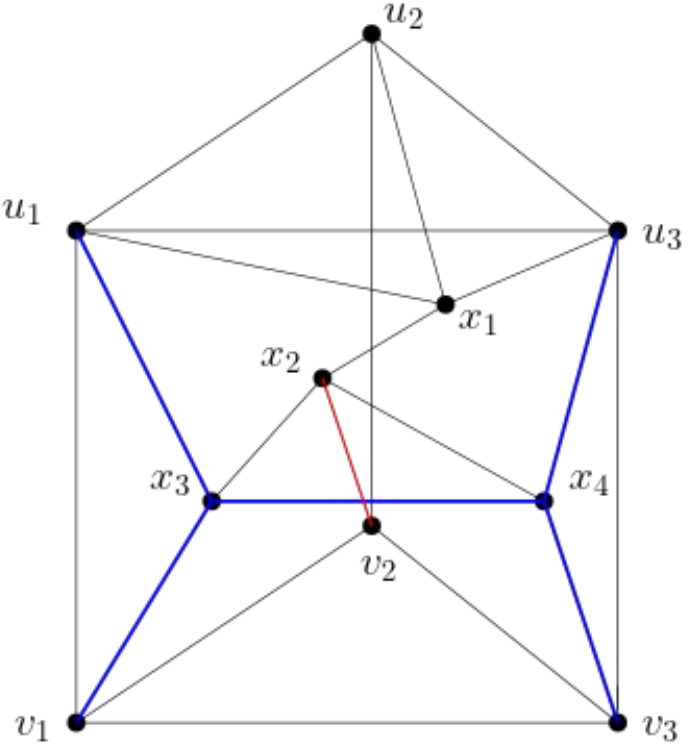}
        		\captionof{figure}{$G_3$, ``rotated" $F'$, case 2}
        		\label{fig:G3_2_2}
        	\end{minipage}
        \end{figure}
        

        It is not hard to see that $x_2 \notin F'$. Otherwise, the subgraph $G[x_1 \cup x_2 \cup U]$ is indecomposable, and the other vertex in $F'$ is a simple extension and the neighbors of $x_2$ are all simple extensions, making $P$ indecomposable, see Figure~\ref{fig:G3_4}.
	 	\begin{figure}[H]
	 	\centering
	 	\includegraphics[width=0.35\textwidth]{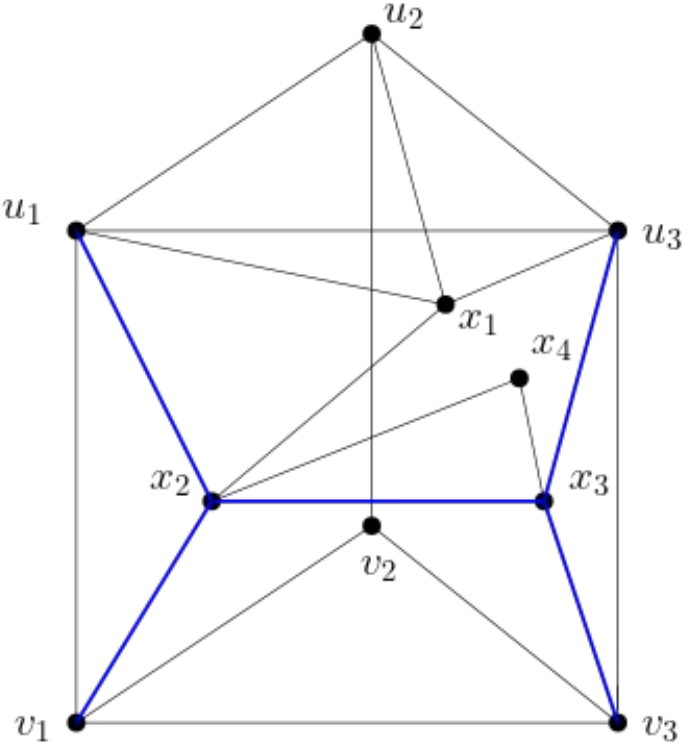}
	 	\caption{$G_3$, ``rotated" $F'$}
	 	\label{fig:G3_4}
	     \end{figure}
     
        Finally, we can assume that every other facet $F'$ that intersects $F$ at a quadrilateral face is not ``rotated". We first show that no triangle edge can be the edge that is associated with such $F'$. If not, say $[x_{i2}, x_{i3}]$ is the edge associated with $F'$, where $x_{i2}$ has at least two neighbors in $U$ and $x_{i3}$ has at least two neighbors in $V$; see Figure~\ref{fig:G3_5}. The third vertex $x_{i4}$ in the triangle must be adjacent to either a vertex in $U$, or a vertex in $V$. Without loss of generality, we may assume $x_{i4}$ is adjacent to a vertex in $V$. So, the subgraph $G[x_{i3} \cup x_{i4} \cup U]$ is indecomposable. Then $x_{i2}$ is a simple extension of this indecomposable subgraph. It follows that the neighbors of $x_{i2}$ in $U$ are simple extensions of $x_{i2} \cup x_{i3} \cup x_{i4} \cup V$. It also follows that $P$ would be indecomposable, a contradiction. Figure~\ref{fig:G3_5} gives an illustration where $[x_2, x_3]$ is the edge associated with $F'$.        
       	\begin{figure}[H]
        	\centering
        	\includegraphics[width=0.45\textwidth]{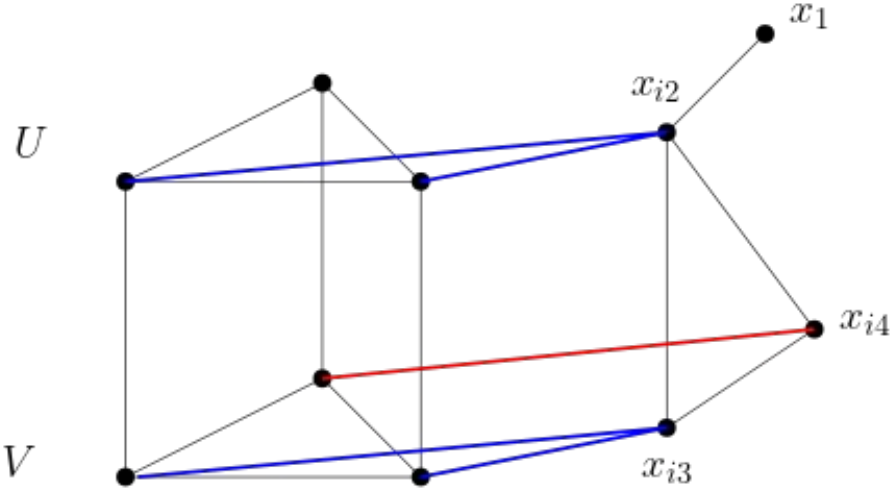}
        	\caption{$G_3$, not ``rotated" $F'$}
        	\label{fig:G3_5}
        \end{figure}
        Now we can confirm that $[x_1, x_2]$ is the edge that is associated with every such $F'$ that intersects $F$ at a quadrilateral face and that the lines containing $[x_1, x_2]$ are projectively parallel to the ``side" edges of $F$. The vertices $x_1, x_2$ cannot be adjacent to both vertices in $U$ and vertices in $V$. And by the argument about skew lines \cref{lit}(6), if $x_2$ is not adjacent to vertices in $U$ (or vertices in $V$, respectively), $x_3, x_4$ are not adjacent to vertices in $U$ (or vertices in $V$, respectively). By \cref{decompfunc-linesegment}, $P$ has a line segment for a summand.
        
	    \textcolor{blue}{$G_4$.} 
	    
	    Label the vertex with degree 3 in $G_4$ as $x_4$ and label the other three vertices as $x_1, x_2, x_3$. As before, we can see that $x_1, x_2, x_3$ are adjacent to vertices either only in $U$ or only in $V$. So every other facet $F'$ that intersects $F$ at a quadrilateral face is not ``rotated". So, we may assume that  $x_1$ is adjacent to $U$, that $x_4$ is adjacent to at least two vertices in $V$, and $[x_1, x_4]$ is projectively parallel to the ``side" edges of $F$; see Figure~\ref{fig:G4_1}. By the argument about skew lines \cref{lit}(6), we can see that $x_4$ is not adjacent to vertices in $U$. It remains to show that $x_2, x_3$ are adjacent to vertices in $V$. Suppose that $x_2$ is not adjacent to any vertices in $U$. Then $x_2$ is adjacent to the 3 vertices in $U$, implying that the subgraph $G[x_1 \cup x_2 \cup U]$ is indecomposable. And the line containing $[x_2, x_4]$ would be skew to a ``side" edge of $F$, making $P$ indecomposable. So now we have a partition of vertices of the polytope, $T_1=\{ x_1, U \}, T_2=\{ x_2, x_3, x_4, V \}$. Additionally, there are 4 edges connecting $T_1, T_2$, whose containing lines are projectively parallel. By \cref{decompfunc-linesegment}, $P$ has a line segment for a summand.
        
	    \begin{figure}[H]
		    \centering
     	    \includegraphics[width=0.45\textwidth]{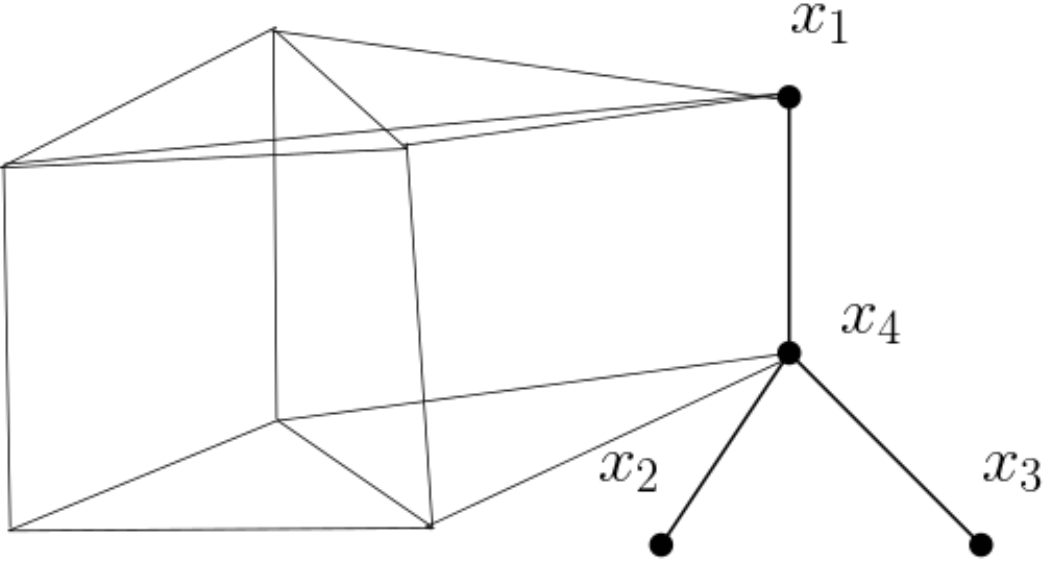}
	    	\caption{$G_4$}
		    \label{fig:G4_1}
     	\end{figure}
     
		\blue{$G_5, G_6$.} 
		
		Since $G_5$ is indecomposable, $F$ has Shephard's property. So $f_1(P)=f_1(F)+\#$edges in between $F$ and $G_5 + \#$edges in $G_5=9+6+5=20$. Similarly, since $G_6$ is indecomposable, $F$ has Shephard's property, so $f_1(P)=9+6+6=21$. Again, by \cite[Theorem 6.1]{PUY2018Excess}, only $A(4), B(4), C(4)$ and $\Sigma(4)$ are possible. However, they each have a decomposable facet which is not a prism, contradicting our assumption that $F$ is a decomposable facet with the maximum number of vertices. So, $G_5, G_6$ are not possible.
	\end{proof}
\end{lemma}

Finally we are ready to present the theorem in dimension 4.

\begin{theorem}\label{thm:dim4}
	Let $P$ be a decomposable $4$-polytope with $10$ vertices. Then $P$ is either $\Sigma(4), C(4)$, or $P$ has a line segment for a summand where a cross section of $P$ is either a $3$-simplex, a bi-pyramid over a triangle, or a square pyramid. In particular, in the case that $P$ has a line a segment for a summand, there is a corresponding partition on the vertices $V(P)=V_1\cup V_2$ with edges between them projectively parallel, where $G[V_1]$, $G[V_2]$ are indecomposable subgraphs.    
	\begin{proof}
		Let $F$ be a decomposable facet of $P$ with the maximum number of vertices.
		
		Case 1: $F$ has 8 vertices. By \cref{dim4:linesegment}, $P$ has a line segment for a summand. Moreover, there is a decomposing function $f$ that maps vertices of $F$ to vertices $u, v$ which are outside $F$. Notice that the induced subgraph of $f^{-1}(u)$ otherricted on $F$ is connected, so the induced subgraph of $f^{-1}(u)$ forms a chain of triangles, therefore $G[V_1]$ is indecomposable. The same argument works for $G[V_2]$.

		Case 2: $F$ has 7 vertices. By \cite[Theorem 5.9]{PUY2018Excess}, $F$ can be a pentasm, a 1-capped prism, or $\Sigma(3)$. If $F$ is $\Sigma(3)$, by \cref{dim4:sigma3}, $P$ is either $\Sigma(4)$ or $C(4)$, and if none of the $F$ is $\Sigma(3)$, by \cref{dim4:nonesigma3}, $P$ has a line segment for a summand. In the scenario that $P$ has a line segment for a summand, a decomposing function gives the partition $V(P)=V_1\cup V_2$, and facet $F$ is either a capped prism or a pentasm, where $G[V_1(F)], G[V_2(F)]$ are obviously indecomposable. Since each vertex outside $F$ has degree at most 2 in the induced subgraph by the vertices outside $F$, each of these vertices (say $u\in V_1$) is adjacent to at least two vertices in $V_1(F)$, in other words, $u$ is a simple extension of $G[V_1(F)]$, making $G[V_1(F)] \cup u$ indecomposable. Hence $G[V_1], G[V_2]$ are indecomposable subgraphs.

		Case 3: $F$ has 6 vertices. All the decomposable facets of $P$ are prisms. By \cref{dim4:allprisms}, $P$ has a line segment for a summand. Denote by $V_1(F)\subset V_1$, $V_2(F)\subset V_2$. By the proof of \cref{dim4:allprisms}, there are four possible induced subgraphs of vertices outside $F$. In $G_1, G_2, G_3$ and $G_4$, three of the vertices have degree at most 2 in the induced subgraphs, so they are adjacent to at least two vertices in the corresponding $V_i(F), i=1,2$, hence they are simple extensions of an indecomposable subgraph. In $G_3, G_4$, there is a vertex (say $v\in V_1$) with degree 1 in the subgraph, but it is adjacent to at most one vertex in $V_2$, therefore, we know that $v$ is a simple extension to $G[V_1\setminus v]$. Hence $G[V_1], G[V_2]$ are indecomposable.

	\end{proof}
\end{theorem}

The discussion in dimension 5 is similar. 

\begin{lemma}\label{dim5:Delta22}
	Let $P$ be a decomposable $5$-polytope with $12$ vertices. Suppose that $F$ is a decomposable facet with the maximum number of vertices and that $F$ is a $\Delta(2,2)$, then $P$ is $\Delta(2,3)$.
	\begin{proof}
		There are exactly 3 vertices $x_1, x_2, x_3$ outside $F$. Suppose $x_1, x_2, x_3$ are pairwise adjacent, then $F$ has Shephard's property (by \cref{lit}(7)), so the number of edges of $P$ is $18+9+3=30$, i.e. $P$ is simple. Then by \cite[Lemma 2.19]{PUY2018Excess}, $P$ is $\Delta(2,3)$. 
		
		Now we may assume $x_1, x_2, x_3$ are not pairwise adjacent. Let $R$ be a prism face of $F$. Consider the other facet $F'$ of $P$ that intersects $F$ at $R$. Note that the facet $F'$ has to be decomposable, otherwise we have an indecomposable subgraph that leaves out at most 2 vertices outside, implying that $P$ would be indecomposable. Since $x_1, x_2, x_3$ are not pairwise adjacent, $F'$ is not $\Delta(2,2)$, so $F'$ is either a capped prism, a pentasm, or a prism. 
		
		If $F'$ is a capped prism, then there are 3 vertices outside $F'$, and they are pairwise adjacent. By \cref{decomposable:triangleoutside}, $F'$ should have a triangle for a summand, which is not possible, by \cref{lem:cappedpentasm_notrianglesummand}.  
		
		If $F'$ is a pentasm, then again there are 3 vertices outside $F'$, and they are pairwise adjacent. Again, by \cref{decomposable:triangleoutside}, $F'$ should have a triangle for a summand, which is not possible by \cref{lem:cappedpentasm_notrianglesummand}. 
		
		The remaining case is that all such $F'$ are prisms. Denote $R_1, R_2, R_3$ the prism ridges of $F$ as in Figure~\ref{fig:dim5_delta22}. We may assume $[x_1, x_2]$ is the edge that is associated with the ridge $R_1$; see Figure~\ref{fig:dim5_delta22}. Then $[x_2, x_3]$ has to be the edge that is associated with the ridge $R_2$, since otherwise, $x_1$ is adjacent to two triangles in $F$, by simple extensions, we will have an indecomposable subgraph that leaves out at most 2 vertices outside. But then, concerning $R_3$, by a similar argument, the edge $[x_1, x_2]$ cannot be associated with $R_3$, the edge $[x_2, x_3]$ cannot be associated with $R_3$, and the vertices $x_1, x_3$ cannot be in $F'$ which are associated with $R_3$, as the two vertices outside a facet have to be adjacent. So, this case does not arise. 
	    \begin{figure}[H]
			\centering
			\includegraphics[width=0.5\textwidth]{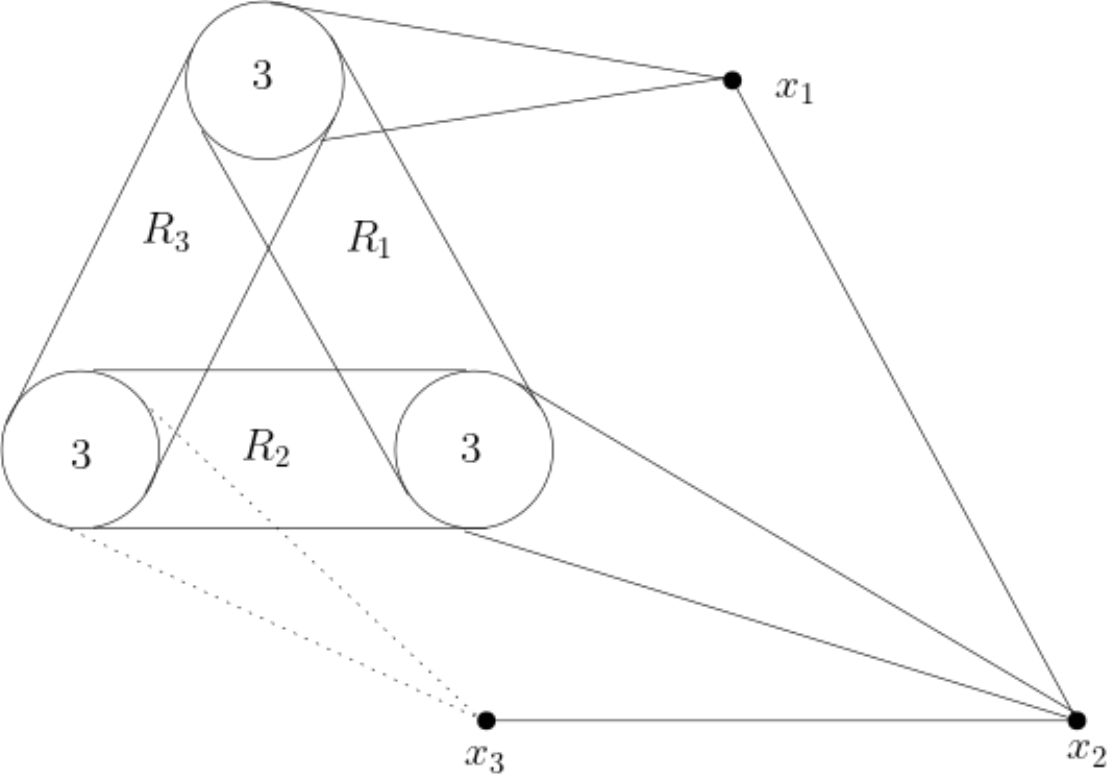}
			\caption{$\Delta(2,2)$ as a facet}
			\label{fig:dim5_delta22}
		\end{figure}
	\end{proof}
\end{lemma}

\begin{lemma}\label{dim5:noneDelta22}
	Let $P$ be a decomposable $5$-polytope with $12$ vertices, and suppose that  $9$ is the maximum number of vertices of any decomposable facet. If no such  facet is $\Delta(2,2)$, then $P$ has a line segment for a summand.
	\begin{proof}
        Let $F$ be such a maximal decomposable facet. Then $F$ is either a capped prism or a pentasm \cite[Theorem 5.9]{PUY2018Excess}. Denote by $U, V$ the partition of vertices of $F$ such that the only edges between them are the projectively parallel ``side" edges of $F$. 
		
		If the three vertices outside $F$ are pairwise adjacent, by \cref{decomposable:triangleoutside}, $F$ has to have a triangle for a summand, which is not possible by \cref{lem:cappedpentasm_notrianglesummand}.
		
		Label the three vertices outside $F$ as $x_1, x_2, x_3$ where $x_1$ is not adjacent to $x_3$. Every other facet $F'$ that intersects $F$ at a ``side" ridge is decomposable, otherwise we would have an indecomposable subgraph that leaves out at most 2 vertices, making $P$ indecomposable. By the assumption, such $F'$ can only be capped prisms, pentasms, or prisms. In any of these cases, there will be an edge outside $F$ that is projectively parallel to the ``side" edges of $F$. Assume without loss of generality that $[x_1, x_2]$ is associated with such an $F'$ and that $x_1$ is adjacent to vertices in $U$, $x_2$ is adjacent to vertices in $V$. It is not hard to check that $x_3$ is only adjacent to vertices in $V$, and $x_1$ is not adjacent any vertex in $V$, and $x_2$ is not adjacent to any vertex in $U$. This gives the required partition of the vertices of $P$, and by \cref{decompfunc-linesegment}, $P$ has a line segment for a summand.		
	\end{proof}
\end{lemma}

\begin{lemma}\label{dim5:prismfacet}
	Let $P$ be a decomposable $5$-polytope with $12$ vertices. If all the decomposable facets of $P$ have $8$ vertices, then $P$ has a line segment for a summand.
	\begin{proof}
		Let $F$ be a simplicial prism which is a facet of $P$, and denote $U, V$ the set of vertices of the simplex facets of $F$, respectively. There are 4 vertices that are outside $F$. Again, there are only 6 possible subgraphs induced by these 4 vertices; see Figure~\ref{fig:prism}. We can see that for the prism 3-faces in $F$, all the other facets $F'$ that intersect $F$ at a prism ridge are decomposable (and hence are all prisms); otherwise we can find an indecomposable subgraph with at most 3 vertices outside, and $P$ would be indecomposable.	
		
		By the structure of 4-prisms, we can now be certain that $F'$ is not ``rotated". It suffices to give a partition of vertices of $P$ such that the edges connecting the two parts are projectively parallel. 
		
		Claim: If $x_i$ (any vertex outside $F$) is adjacent to at least 3 vertices in $F$, then $x_i$ is either only adjacent to vertices in $U$ or only to vertices in $V$. 
		
		Suppose the contrary. Then we can find a 5-cycle, with vertices $x_i$, the two neighbors of $x_i$ in $F$ which are not adjacent, and two other vertices of $F$ which are adjacent; see Figure~\ref{fig:dim5_claim1} for illustration. The affine dimension of this 5-cycle is 4. So, by \cite[Proposition 2]{Przeslawski-Yost-2008}, this 5-cycle is indecomposable. Since the 5-cycle $C$ shares a common edge with the subgraph induced by $U$ and a common edge in the subgraph induced by $V$, which are both indecomposable, the subgraph $G[U \cup V] \cup G[C]$ is indecomposable, so we have an indecomposable subgraph that leaves out at most 3 vertices, making $P$ indecomposable, a contradiction.
	   \begin{figure}[H]
			\centering
			\includegraphics[width=0.35\textwidth]{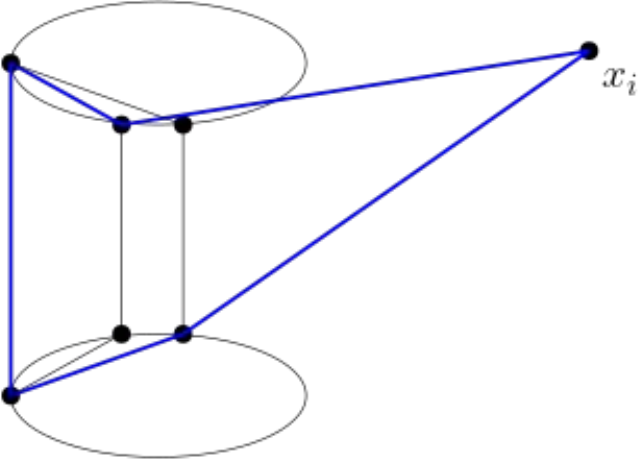}
			\caption{a 5-cycle of dimension 4}
			\label{fig:dim5_claim1}
		\end{figure}
	
	    Now we discuss the subgraphs induced by the 4 vertices outside $F$. 
	    
	    In $G_1$, we may assume $[x_1, x_2]$ is an edge associated with the other prism facet $F'$. To show that $P$ has a line segment for a summand, it suffices to show that $x_4$ is only adjacent to vertices in $U$ and $x_3$ is only adjacent to vertices in $V$. By the claim, we just need to exclude the situation where $x_4$ is only adjacent to vertices in $V$. 
	    Now suppose that $x_4$ is only adjacent to vertices in $V$. Since $G[x_1 \cup U]$ is indecomposable, $G[x_2 \cup x_4 \cup V]$ is indecomposable. However, since the line containing $[x_1, x_4]$ is skew to some of the lines containing the ``side" edges of $F$, we have an indecomposable subgraph that leaves only 1 vertex outside, making $P$ indecomposable. A similar argument applies for $x_3$. Thus, we have a partition of vertices of $P$, $T_1=\{x_1, x_4, U\}$, $T_2=\{x_2, x_3, V\}$. As both of the subgraphs induced by $T_1, T_2$ are indecomposable, the edges connecting them are projectively parallel. So, by \cref{decompfunc-linesegment}, $P$ has a line segment for a summand. 
	    \begin{figure}[H]
	    	\centering
	    	\includegraphics[width=0.4\textwidth]{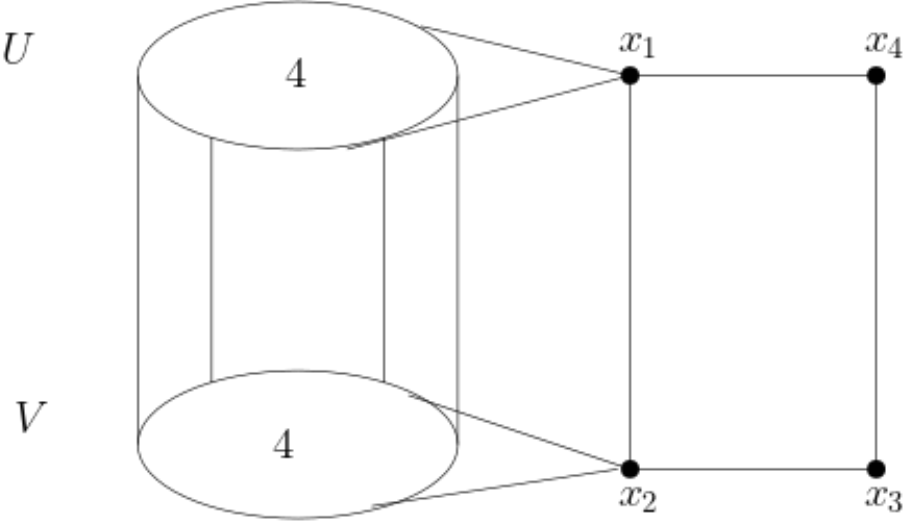}
	    	\caption{$G_1$, dimension 5}
	    	\label{fig:dim5_G1}
	    \end{figure}
    
	    In $G_2$, again, we know that there is an edge, without loss of generality say $[x_1, x_2]$, outside $F$ that is associated with the other facet $F'$. Assume without loss of generality that $x_1$ is adjacent to vertices in $U$ and $x_2$ is adjacent to vertices in $V$. By the discussion related to $G_1$, since $x_3$ is adjacent to $x_1$, the vertex $x_3$ is adjacent to vertices in $U$ but not vertices in $V$. By the claim, the vertex $x_4$ is adjacent to either vertices in $U$ but not vertices in $V$ or is adjacent to vertices in $V$ but not $U$. In either of these cases, we can have a partition of the vertices of $P$ where $G[U]$ and $G[V]$ are indecomposable. Additionally, the edges in between the partition have to be projectively parallel, reasoned as before. By \cref{decompfunc-linesegment}, $P$ has a line segment for a summand. 
	    
	    In $G_3$, denote the three vertices of the triangle as $x_{i1}, x_{i2}, x_{i3}$. We first show that an edge of the triangle cannot be the edge that is associated with $F'$ where $F'$ intersects $F$ at a ``side'' ridge. Suppose not, say $[x_{i1}, x_{i2}]$ is the edge that is associated with $F'$, where $x_{i1}$ is adjacent to vertices in $U$, and $x_{i2}$ is adjacent to vertices in $V$. Since $x_{i3}$ is adjacent to some vertices in $U$ or (say) some vertices in $V$, the subgraph $G[x_{i2} \cup x_{i3} \cup V]$ is indecomposable and the line containing $[x_{i1}, x_{i3}]$ is skew to the ``side" edge of $F$ that connects $G[x_{i2} \cup x_{i3} \cup V]$ and $G[x_{i1} \cup U]$ and so, $P$ would be indecomposable, a contradiction; see Figure~\ref{fig:dim5_G3}. So, $[x_{i1}, x_{i4}]$ is the unique edge that is associated with all the other facets $F'$. And it is not hard to see that the adjacency between $x_{i2}, x_{i3}$ and $U, V$ are consistent with the adjacency of $x_{i1}$, in the sense of either only vertices of $U$, or only vertices of $V$. By \cref{decompfunc-linesegment}, $P$ has a line segment for a summand.
	    \begin{figure}[H]
	    	\centering
	    	\includegraphics[width=0.4\textwidth]{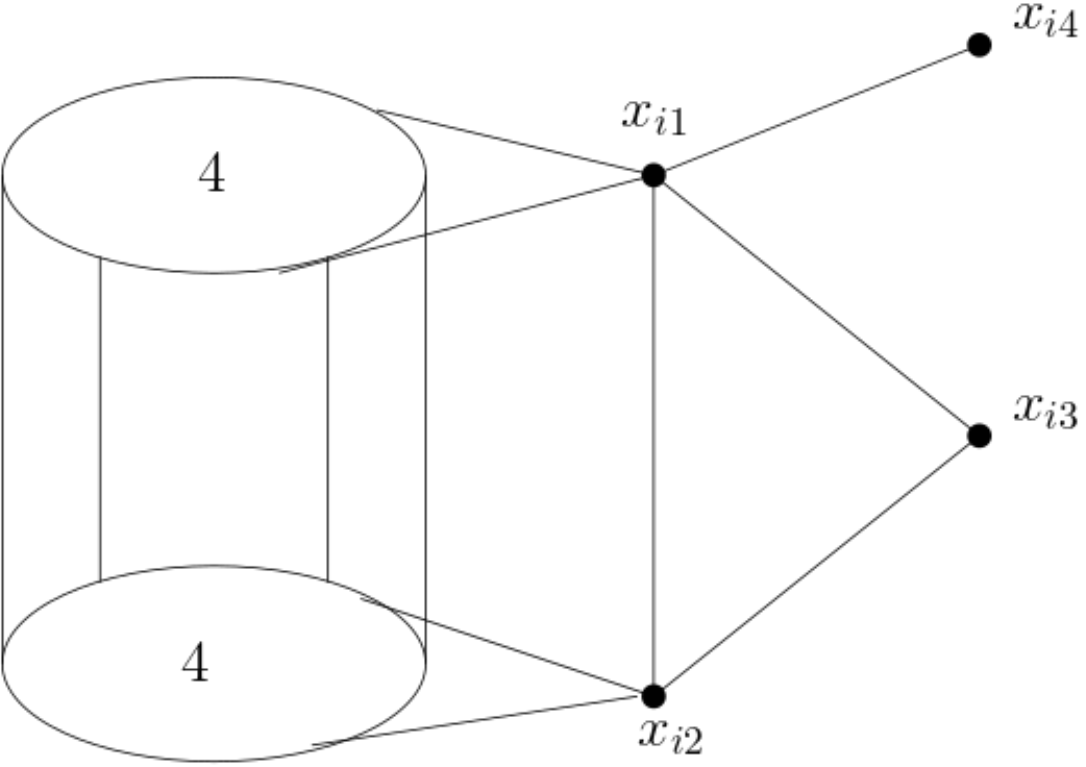}
	    	\caption{$G_3$, dimension 5}
	    	\label{fig:dim5_G3}
	    \end{figure}
    
	    In $G_4$, the three vertices  with degree 1 are in a symmetric situation, so label them as $x_2, x_3, x_4$. And we may assume that $[x_1, x_2]$ is an edge that is associated with $F'$. By the discussion related to $G_1$, the adjacency of $x_3, x_4$ to $U, V$ is consistent with the adjacency of $x_1$ to $U, V$. So, this gives a partition of the vertices of $P$ and the edges connecting the partition sets are projectively parallel. By \cref{decompfunc-linesegment}, $P$ has a line segment for a summand. 
	    
	    In $G_5, G_6$. Since $G_5, G_6$ are indecomposable, $F$ has Shephard's property. So, if it were $G_5$, the number of edges of $P$ would be 29, so $G_5$ is not possible; if it were $G_6$, the number of edges of $P$ would be 30, i.e. $P$ would be simple, so, $P$ would be $\Delta(2,3)$. However, $\Delta(2,3)$ has a decomposable facet with 9 vertices, contradicting our assumption, so $G_6$ is not possible either.
	\end{proof}
\end{lemma}


\begin{theorem}\label{thm:dim5}
	Let $P$ be a decomposable $5$-polytope with $12$ vertices. Then $P$ is either $\Delta(2,3)$ or has a line segment for a summand whose cross section is either a $4$-simplex or a $4$-polytope with $6$ vertices.
	\begin{proof}
     	Let $F$ be a decomposable facet of $P$ with the maximum number of vertices.
     	
		Case 1: $F$ has 10 vertices. By \cref{2verticesoutside}, $P$ has a line segment for a summand. 
		
		Case 2: $F$ has 9 vertices. By \cite{PUY2018Excess}, $F$ is either a pentasm, a capped prism, or $\Delta(2,2)$. If $F$ is $\Delta(2,2)$, by \cref{dim5:Delta22}, $P$ is $\Delta(2,3)$; if none of the $F$ is $\Delta(2,2)$, by \cref{dim5:noneDelta22}, $P$ has a line segment for a summand. 
		
		Case 3: $F$ has 8 vertices. All the decomposable facets are simplicial prisms. By \cref{dim5:prismfacet}, $P$ has a line segment for a summand.
		
		When $P$ has a line segment for a summand, by \cref{cross-section}, we can take a cross section of $P$, and the result follows. 
	\end{proof}
\end{theorem}

\begin{corollary}\label{corollary:2dplus2}
	Let $P$ be a decomposable $d$-polytope with $2d+2$ vertices and have a line segment for a summand, $d\ge 4$. Then there is a corresponding partition of the vertices of $V(P)=V_1 \cup V_2$ with edges between $V_1$ and $V_2$ are projectively parallel and the subgraphs of $G[V_1], G[V_2]$ are both indecomposable.
	\begin{proof}
    It suffices to show that the subgraphs $G[V_1]$ and $ G[V_2]$ are indecomposable. We proceed by induction on dimension. The induction base is dimension 4, which is proved by \cref{thm:dim4}. We may let $V_1, V_2$ be a partition of vertices of $P$ such that edges between $V_1$ and $V_2$ are ``side" edges. Choose a ``side" facet $F$ (i.e. a facet contains some of the ``side" edges). Then there are at most 4 vertices outside $F$. Denote by $V_1(F) \subset V_1$ and $V_2(F) \subset V_2$ the corresponding partition of vertices in $F$. Choose a vertex $v\in V(P) \setminus V(F)$. The vertex $v$ has at least 2 neighbors in $V_i(F)$, so $v$ is a simple extension of the subgraph induced by $V_i$. Thus, the subgraph $G[v \cup V_i]$ is indecomposable. The result now follows.
	\end{proof}
\end{corollary}

Finally, we are ready to give the results for all dimensions higher than 4. As remarked in the introduction, $F287$ shows that the next two results fail in dimension 3.

\begin{theorem}\label{thm:decomposable2dplus2}
	Let $P$ be a decomposable $d$-polytope with $2d+2$ vertices, $d\ge 4$. Then $P$ is either $C(4), \Sigma(4)$, $\Delta(2,3)$ or has a line segment for a summand whose cross section is a $(d-1)$-polytope with at most $d+1$ vertices.	
	\begin{proof}
		Let $F$ be a decomposable facet with maximum number of vertices. Then $F$ has $2d-2, 2d-1$ or $2d$ vertices. It suffices to prove the statement when $d\ge 6$. 
		
		Case 1: $F$ has $2d=2(d-1)+2$ vertices. By \cref{2verticesoutside}, $P$ has a line segment for a summand. 
		
		Case 2: $F$ has $2d-1=2(d-1)+1$ vertices. By \cite[Theorem 5.9]{PUY2018Excess}, $F$ is either a pentasm or a capped prism. By the decomposability of $P$, the other facet $F'$ that intersects $F$ at a pentasm or capped prism ridge is also decomposable, so $F'$ is either a pentasm, a capped prism or a prism. Denote a partition of the vertices of $F$ by $U,V$ as before, so that all edges between vertices in $U$ and vertices in $V$ are projectively parallel. By the structure of $F'$, there is an edge $[x_1, x_2]$ outside $F$ that is projectively parallel to the ``side" edges of $F$ and the subgraphs $G[x_1 \cup U]$, $G[x_2 \cup V]$ are indecomposable. Denote by $x_3$ the other vertex outside $F$. Since $d\ge 6$, the vertex $x_3$ is adjacent to at least 4 vertices in $F$, we can see that $x_3$ is adjacent to either (say) only vertices in $U$, or only vertices in $V$. So, we have a partition of vertices of $P$, $T_1=\{ x_1, x_3, U  \}$, $T_2= \{ x_2, V  \}$. The subgraphs $G[T_1], G[T_2]$ are indecomposable and so, all the edges in between them are projectively parallel. By \cref{decompfunc-linesegment}, $P$ has a line segment for a summand. 
		
		Case 3: $F$ has $2d-2=2(d-1)$ vertices, i.e. all the decomposable facets are prisms. Again, there are six possible subgraphs induced by the 4 vertices outside $F$; see Figure~\ref{fig:prism}. By the claim in \cref{dim5:prismfacet}, one can see that in $G_1, G_2, G_3, G_4$, all the vertices are only adjacent to one side of the vertices of $F$, which gives the partition of the vertices of $P$ and illustrates the indecomposability of the induced subgraph of each partition. So the edges in between each partition are all projectively parallel, implying that $P$ has a line segment for a summand (by \cref{decompfunc-linesegment}). As for $G_5$ and $ G_6$, since $G_5, G_6$ are indecomposable, $F$ would have Shephard's property. Thus, the number of edges of $P$ is $2\binom{d-1}{2}+(d-1)+2(d-1)+5=d^2+4$ or $2\binom{d-1}{2}+(d-1)+2(d-1)+6=d^2+5$. A $d$-polytope with $2d+2$ vertices should have at least $d^2+d >d^2+5$ edges when $d\ge 6$. Thus $G_5, G_6$ are not possible. 
		
		When we know that $P$ has a line segment for a summand, by \cref{cross-section}, we can take a suitable cross section $H \cap P'$ where $P'$ is projectively equivalent to $P$ where $H \cap P'$ necessarily has either $d$ vertices or $d+1$ vertices.		
	\end{proof}
\end{theorem}

\begin{remark}
	The polytope $\Delta(2,3)$ does not appear as a facet of a decomposable $d$-polytope with $2d+2$ vertices. If it occurred, it would be a facet of a 6-polytope with 14 vertices, then $\Delta(2,3)$ would have a line segment for a summand, a contradiction.
\end{remark}

By \cite[Theorem~4]{WJDYconditional}, the structure of decomposable $d$-polytopes with $2d+2$ vertices gives the combinatorial decomposability.

\begin{corollary}\label{cor:2dplus2decomposable}
	There is no conditionally decomposable $d$-polytope with $2d+2$ vertices, for all $d\ge 4$.
\end{corollary}


\section{Decomposable $d$-polytopes with $2d+3$ vertices, $d\ge 5$}

\begin{lemma}\label{dim5:Sigma4}
	Let $P$ be a decomposable $5$-polytope with $13$ vertices and let $F$ be a decomposable facet with the maximum number of vertices. If $F$ is $\Sigma(4)$, then $P$ is $\Sigma(5)$.
	\begin{proof}
		Consider a ridge $R$ which is $\Sigma(3)$ in $F$. It is routine to check and see that the other facet $F'$ that intersects $F$ at $R$ is decomposable. So $F'$ is a decomposable 4-polytope with 8 ($2d$), 9 ($2d+1$), or 10 ($2d+2$) vertices that has  $\Sigma(3)$ as a facet. So, $F'$ is either $\Sigma(4)$ or $C(4)$. Observe first that, if one of such $F'$ is $\Sigma(4)$, then all such $F'$ are $\Sigma(4)$; if one of such $F'$ is $C(4)$, then all such $F'$ are $C(4)$. This is because the 3 vertices outside a $\Sigma(3)$ facet in $C(4)$ are pairwise adjacent and the 3 vertices outside a $\Sigma(3)$ facet in $\Sigma(4)$ are not pairwise adjacent, and one facet $F'$ determines the adjacency between these 3 vertices. 
		
		Denote the unique nonsimple vertex of $F$ (=$\Sigma(4)$) as $O$, denote the vertex set of the triangular faces of $F$ as $V_1, V_2, V_3$, and denote the 3 vertices outside $F$ as $a,b,c$. 
		
		Suppose $F'$ is a $\Sigma(4)$; see Figure~\ref{fig:dim5_Sigma4_Sigma5}. W.l.o.g. we may assume that $a$ is not adjacent to $c$, and that the vertex $a$ is adjacent to $O$ and to vertices in $V_1$, and hence adjacent to every vertex in $V_1$; the vertex $c$ is adjacent to $O$ and vertices in $V_3$, and hence adjacent to every vertex in $V_3$; the vertex $b$ is adjacent to vertices in $V_2$, and hence adjacent to every vertex in $V_2$. This together with the graph of $F$ gives the entire graph of $P$, where the polytope $P$ has excess degree $3$ and $O$ is the only nonsimple vertex in $P$. Thus, by \cite[Theorem 4.9]{PUY2020} (polytopes with excess degree at most $d-1$ are determined by their edge graph), $P$ is $\Sigma(5)$.
	    \begin{figure}[H]
			\centering
			\includegraphics[width=0.45\textwidth]{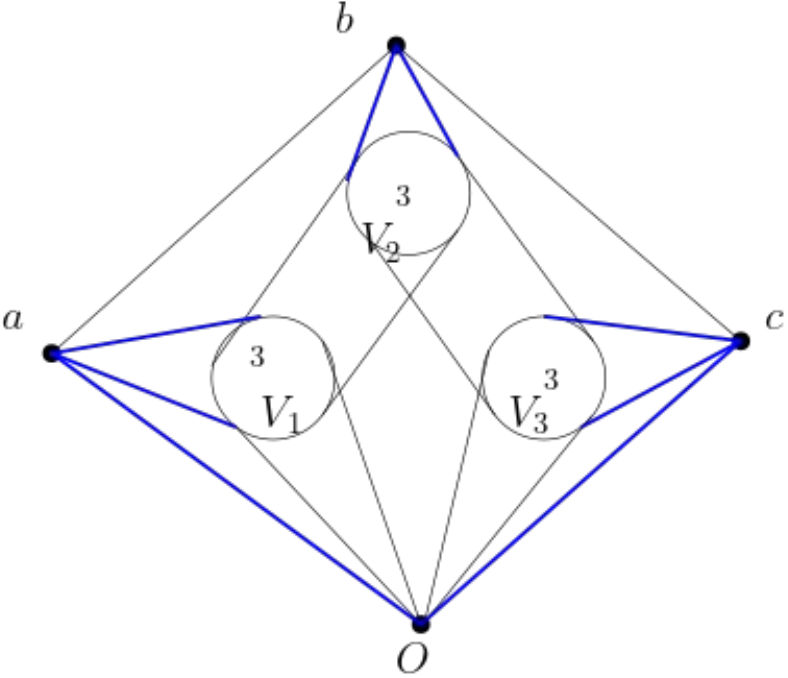}
			\caption{$F'$ is $\Sigma(4)$}
			\label{fig:dim5_Sigma4_Sigma5}
		\end{figure}
	
		Suppose $F'$ is a $C(4)$. Then the three vertices $a, b, c$ outside $F$ are pairwise adjacent and $F$ has Shephard's property; see Figure~\ref{fig:dim5_C4_Sigma5}. By \cref{decomposable:triangleoutside}, $F$ has a triangle for a summand and the adjacency of the vertices of $F$ with $a, b, c$ defines a decomposing function $f$ for $F$. Let us have a look at $F'$ (which is $C(4)$). There is a 1-1 correspondence between the triangular face of $\Sigma(3)$ and $\triangle_{a,b,c}$, so vertices in $V_1$ and $O$ are not adjacent to the same vertex in $\triangle_{a,b,c}$. However, a proper decomposing function should map a 3-simplex ($\conv\{O, V_1\}$) to the same vertex of a 2-simplex ($\conv\{a,b,c\}$), this gives a contradiction. So $F'$ cannot be $C(4)$, and this case does not arise.
	    \begin{figure}[H]
			\centering
			\includegraphics[width=0.35\textwidth]{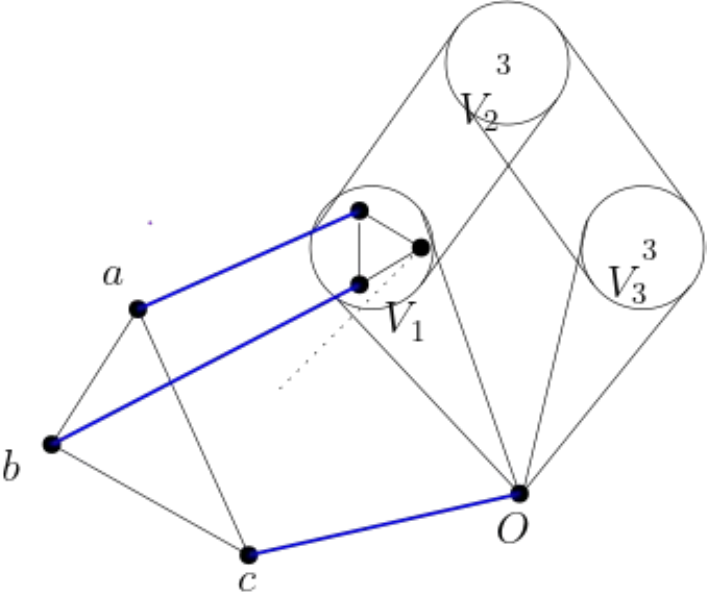}
			\caption{$C(4)$ as a facet}
			\label{fig:dim5_C4_Sigma5}
		\end{figure}		
	\end{proof}
\end{lemma}

\begin{lemma}\label{dim5:C4}
	Let $P$ be a decomposable $5$-polytope with $13$ vertices, and let $F$ be a decomposable facet with the maximum number of vertices. If $F$ is $C(4)$, then $P$ is $C(5)$.
	\begin{proof}
		We apply the similar approach as the previous lemma. Consider a facet $F'$ that intersects $F$ at a $R=\Sigma(3)$. It is routine to check that $F'$ is a decomposable 4-polytope, with 8, 9 or 10 vertices and has $\Sigma(3)$ as a facet. So $F'$ is either $C(4)$ or $\Sigma(4)$. If one such $F'$ is $C(4)$, then they are both $C(4)$; if one such $F'$ is $\Sigma(4)$, then they are both $\Sigma(4)$. 
		
		Suppose $F'$ is $C(4)$. Then this gives the adjacency of $a, b, c$ with almost all vertices in $F$ except for the vertices in a triangle. Since the other facet is also $C(4)$, the adjacency of $a, b, c$ with the vertices of $F$ is clear; see Figure~\ref{fig:dim5_C4_C5_1}. As $P$ has excess degree 3, by \cite[Theorem 4.9]{PUY2020}, the graph of $P$ determines the face lattice of $P$. Thus $P$ is $C(5)$.
	    \begin{figure}[H]
			\centering
			\includegraphics[width=0.4\textwidth]{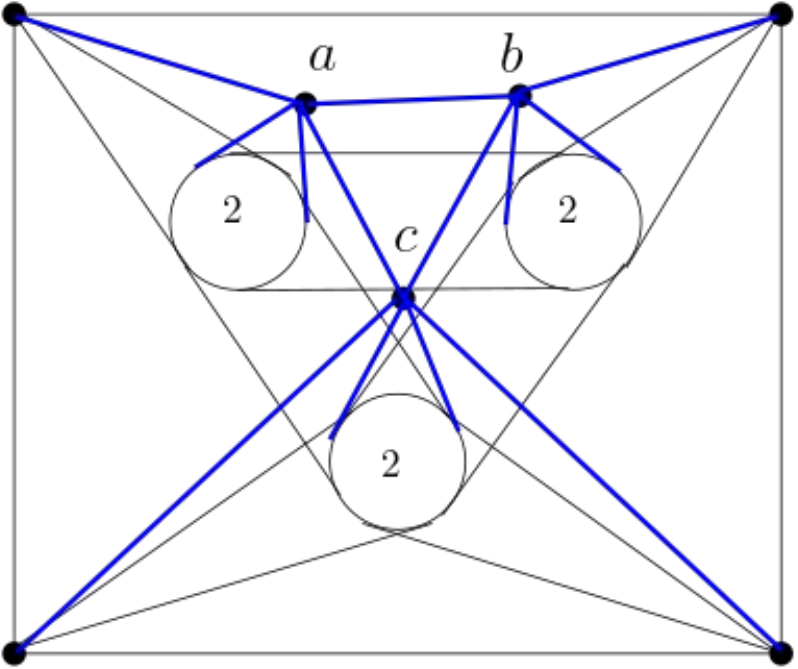}
			\caption{$C(4)$ as a facet}
			\label{fig:dim5_C4_C5_1}
		\end{figure}
	
		Suppose $F'$ is $\Sigma(4)$. By the structure of $\Sigma(4)$, there is a vertex, say $a$, that is adjacent to the bottom triangle of $F$ (refer Figure~\ref{fig:dim5_C4_C5_1}), and there is a vertex $b$ that is adjacent to the upper triangle of $F$ (refer Figure~\ref{fig:dim5_C4_C5_1}), where $a,b$ are not adjacent. Then the edges of the last vertex $c$ is clear, see the blue edges in Figure~\ref{fig:dim5_C4_C5_2}. If we consider the other facet $F''$ which is $\Sigma(4)$, one can see that the adjacency follows the same pattern, i.e. the vertex $b$ is adjacent to the other triangle in $F$, and the vertex $a$ is adjacent to one more vertex in $F$, see the green edges in Figure~\ref{fig:dim5_C4_C5_2}. Note that this gives the adjacency between $a,b,c$ and $F$, hence the entire edge graph of $P$. However, the vertex $c$ has degree 4 in a 5-dimensional polytope, which is impossible. So, this case does not arise.
		\begin{figure}[H]
			\centering
			\includegraphics[width=0.45\textwidth]{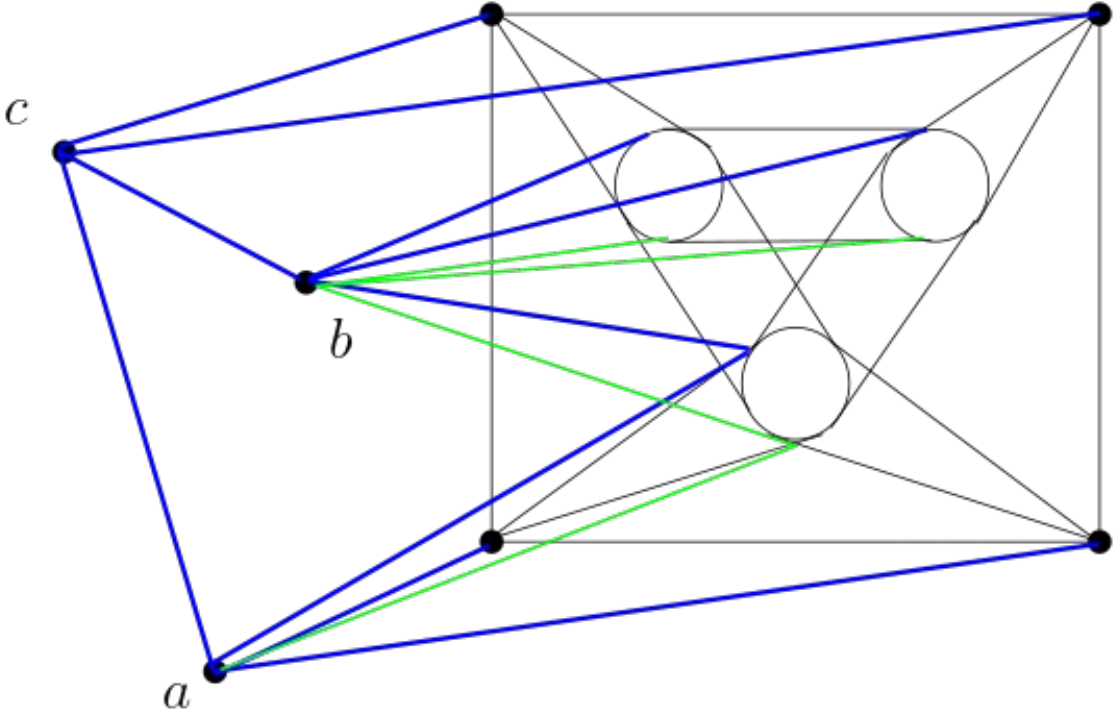}
			\caption{$C(4)$ as a facet}
			\label{fig:dim5_C4_C5_2}
		\end{figure}	
	\end{proof}
\end{lemma}

\begin{lemma}\label{dim5:Delta22_2dplus3}
	Let $P$ be a decomposable $5$-polytope with $13$ vertices, and let $F$ be a decomposable facet with the maximum number of vertices. Then $F$ cannot be $\Delta(2,2)$.
	\begin{proof}
		Suppose $F$ is $\Delta(2,2)$. We will prove the lemma by contradiction. Consider the other facet $F'$ that intersects  $F$ at a prism ridge of $P$. Then it is routine to check that $F'$ is a decomposable facet with 8 or 9 vertices. So, $F'$ could be a prism, pentasm, capped prism, or $\Delta(2,2)$. 
		
		We first show that $F'$ cannot be $\Delta(2,2)$. Suppose the contrary, i.e. $F, F'$ are both $\Delta(2,2)$. Then the subgraph induced by the 4 vertices outside $F$ must contain a triangle and so, it could only be $G_3, G_5$ or $G_6$ from Figure~\ref{fig:prism}. Note that $G_5, G_6$ are not possible, as $F$ would have to have Shephard's property, which would imply that $P$ has 32 or 33 edges, i.e. $P$ has excess degree $-1$ or 1 which is not possible. Denote the vertex that has degree 1 in $G_3$ by $x_4$ and denote the vertex sets of triangular faces of $F$ by $V_1, V_2, V_3$; see Figure~\ref{fig:dim5_Delta22_2dplus3}. Then $x_4$ is adjacent to at least 4 vertices in $F$. By the pigeonhole principle, $x_4$ is adjacent to at least 2 vertices in a $V_i$. So the subgraph $G[x_4 \cup V_i]$ is indecomposable. Since $x_4$ is adjacent to a vertex $v$ in $V_j$ ($j\ne i$), by looking at the edge graph of $\Delta(2,2)$, we can see that $v$ is a simple extension of $x_4 \cup V_i$, implying that the subgraph $G[x_4 \cup v \cup V_i]$ is indecomposable. Hence, the subgraph $G[x_4 \cup V_j \cup V_i]$ is indecomposable. Then, the subgraph $G[x_4 \cup V_1 \cup V_2 \cup V_3]$ is indecomposable, making $P$ indecomposable, a contradiction. 
		\begin{figure}[H]
			\centering
			\includegraphics[width=0.45\textwidth]{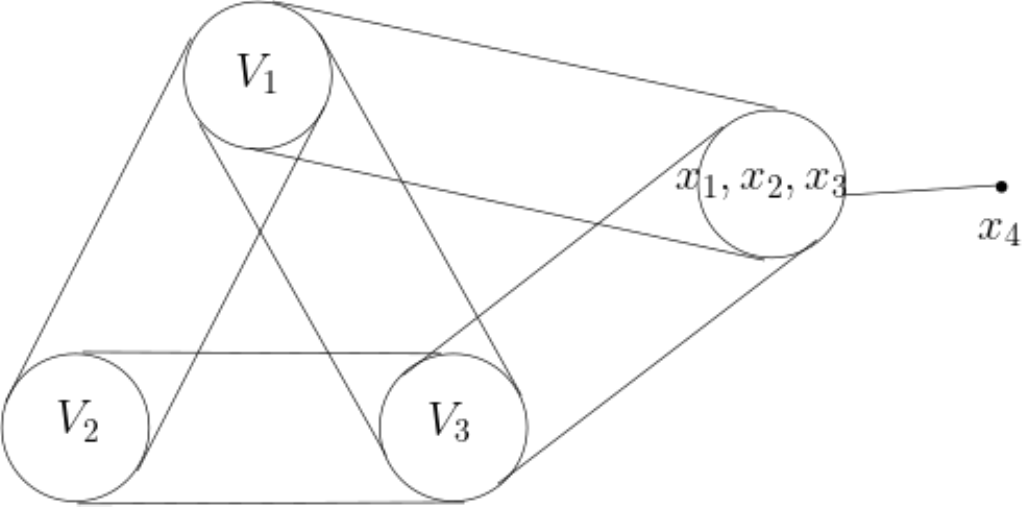}
			\caption{$\Delta(2,2)$ as a facet}
			\label{fig:dim5_Delta22_2dplus3}
		\end{figure}
		The remaining case is to show that $F'$ cannot be a pentasm, a capped prism or a prism. Denote the vertex set of the triangular faces of $F$ as $V_1, V_2, V_3$ and $T_1, T_2, T_3$; see Figure~\ref{fig:dim5_Delta22_2dplus3_2}. Consider a ridge $R_1=\conv\{V_1, V_2\}$. There is an edge (say) $[x_1, x_2]$ outside $F$ that is associated with $R_1$, such that (say) $x_1$ is only adjacent to (at least 2) vertices in $V_1$ and $x_2$ is only adjacent to (at least 2) vertices in $V_2$. So, the subgraphs $G[x_1 \cup V_1]$ and $G[x_2 \cup V_2]$ are indecomposable. It follows that $x_1, x_2$ cannot be adjacent to any vertex in $V_3$, since otherwise their neighbor in $V_3$ would be a simple extension, and then by some simple extensions, we would have an indecomposable subgraph that leaves out at most 3 vertices, making $P$ indecomposable. So, by a similar argument, there is a vertex $x_3$ outside $F$ that is only adjacent to vertices in $V_3$. Consider $R_2=\conv\{T_1, T_2\}$. There is an edge $[x_i, x_j]$ that is associated with $R_2$, where $x_i$ is adjacent to (at least 2) vertices in $T_1$, and $x_j$ is adjacent to (at least 2) vertices in $T_2$. In other words, $x_i$ and $ x_j$ have to be adjacent to some vertices in both $V_i$ and $V_j$ ($i\in \{i,2,3\}$). However, there are in total 4 vertices outside $F$, the vertices $x_1, x_2, x_3$ are only adjacent to vertices in one of the $V_i$, but we need at least 2 vertices that are adjacent to vertices in two of the $V_i$. So, this is not possible.		
		\begin{figure}[H]
			\centering
			\includegraphics[width=0.35\textwidth]{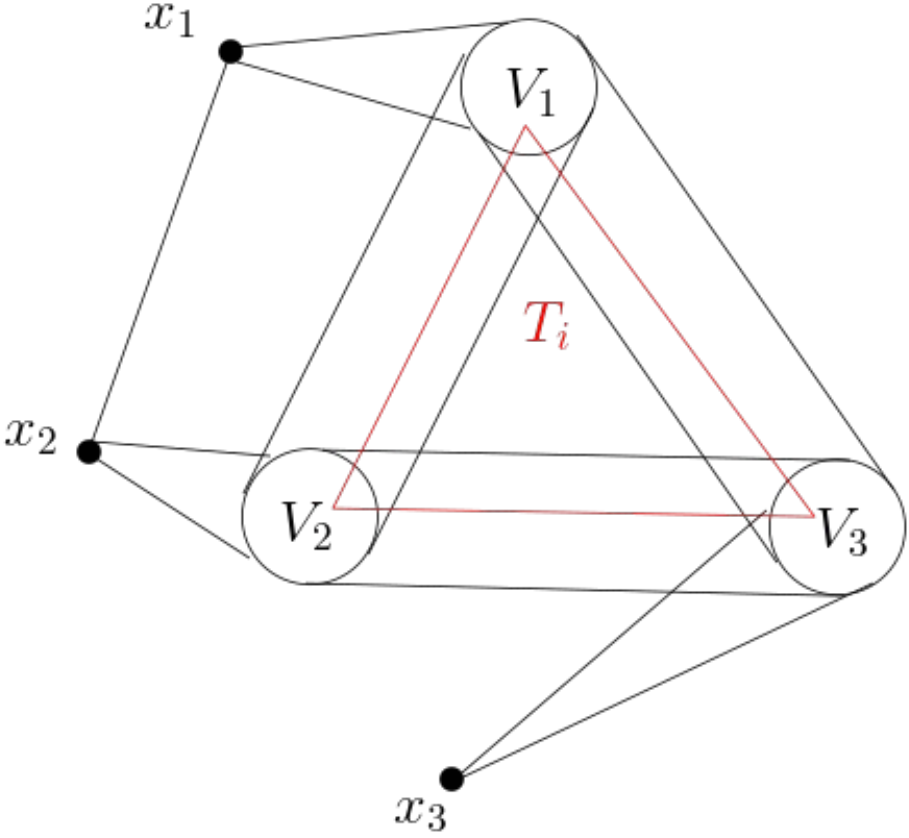}
			\caption{$\Delta(2,2)$ as a facet}
			\label{fig:dim5_Delta22_2dplus3_2}
		\end{figure}	
	\end{proof}
\end{lemma}

\begin{lemma}\label{dim6:delta23}
	Let $P$ be a decomposable $6$-polytope with $15$ vertices and let $F$ be a decomposable facet of $P$ with the maximum number of vertices. If $F$ is $\Delta(2,3)$, then $P$ is $\Delta(2,4)$.
	\begin{proof}
		Consider a ridge $R$ which is $\Delta(2,2)$ in $F$. The other facet $F'$ that intersects $F$ at $R$ must be decomposable; otherwise we have an indecomposable subgraph that leaves out at most 2 vertices (as the remaining three vertices in $F$ are simple extensions) and touches every facet of $P$, which would be a contradiction. So $F'$ is a decomposable 5-polytope with 10 ($2d$), 11 ($2d+1$), or 12 ($2d+2$) vertices and has $\Delta(2,2)$ as a facet. If $F'$ has 10 vertices, $F'$ is a pyramid, so this is not possible. By \cref{lem:Delta22_nolinesegment}, the facet $\Delta(2,2)$ does not have a line segment for a summand, and by \cref{2verticesoutside}, polytopes with a line segment for a summand do not have $\Delta(2,2)$ as a facet, we can conclude by \cref{thm:decomposable2dplus2} that $F'$ is $\Delta(2,3)$. This implies that the 3 vertices outside $F$ (or $F'$ if you want) are pairwise adjacent. Thus, $F$ has Shephard's property. Therefore, the number of edges of $P$ is exactly 45, i.e. $P$ is simple. By the classification of simple polytopes with up to $3d-1$ vertices \cite[Lemma 2.19]{PUY2018Excess}, $P$ is $\Delta(2,4)$.
	\end{proof}
\end{lemma}

\begin{theorem}\label{2dplus3}
	Let $P$ be a decomposable $d$-polytope with $2d+3$ vertices, where $d \ge 5$. Then $P$ is either $\Sigma(5)$, $C(5)$, $\Delta(2,4)$ or has a line segment for a summand, whose cross section is either a $(d-1)$-simplex or a $(d-1)$-polytope with $d+1$ vertices.
	\begin{proof}
		Let $F$ be a decomposable facet of $P$ with the maximum number of vertices. 
		
		Case 1: $F$ has $2d+1$ vertices. By \cref{2verticesoutside}, $P$ has a line segment for a summand. 
		
		Case 2: $F$ has $2d=2(d-1)+2$ vertices. By \cref{thm:decomposable2dplus2}, $F$ is either $C(4), \Sigma(4), \Delta(2,3)$, or has a line segment for a summand. By \cref{dim5:C4}, \cref{dim5:Sigma4}, \cref{dim6:delta23}, if $F$ is $C(4), \Sigma(4)$ or $\Delta(2,3)$, then $P$ is $C(5), \Sigma(5)$ or $\Delta(2,4)$, respectively. Assume now that $F$ has a line segment for a summand, with the partition of vertices $V(F)=V_1(F) \cup V_2(F)$ such that the edges between $V_1(F)$ and $V_2(F)$ are projectively parallel. Consider the other facet $F'$ that intersects $F$ at a ``side" ridge of $P$. It is routine to check that $F'$ is decomposable. We can also assume that $F'$ is not $C(4), \Sigma(4), \Delta(2,3)$. Suppose that $F'$ is a copy of $\Delta(2,2)$ (we will show that this is impossible), then the three vertices $x_1, x_2, x_3$ outside $F$ are pairwise adjacent. Denote the vertices in $F'$ as $x_1, x_2, x_3$, $u_1, u_2, u_3 \in V_1(F)$, $v_1, v_2, v_3 \in V_2(F)$ in a way that  $x_i, u_i, v_i$ are mutually adjacent, for each $i=1,2,3$. Choose a vertex $w\in V_1(F)\setminus \{u_1, u_2, u_3\}$. Since $G[V_1(F)]$ is indecomposable (by \cref{corollary:2dplus2}) and one of the $x_i$ (say $x_1$) has to be adjacent to $w$, we know that $x_1$ is a simple extension of $G[V_1(F)]$. It follows that $x_2, x_3,v_1$ are simple extensions of $G[V_1(F) \cup x_1]$. This would imply that our polytope is indecomposable, a contradiction.
        Now we know that $F'$ also has a line segment for a summand. So, there is an edge $[x_1, x_2]$ outside $F$ that is projectively parallel to the ``side" edges of $F$. We may assume that $x_1$ is adjacent to some (in fact at least 3) vertices in $V_1(F)$. By \cref{corollary:2dplus2}, the subgraph induced by $V_1(F)$ is indecomposable, so the subgraph $G[x_1 \cup V_1(F)]$ is indecomposable. Similarly, the subgraph $G[x_2 \cup V_2(F)]$ is indecomposable. Assume without loss of generality, that the remaining vertex $x_3\notin F$ is adjacent to $x_1$. Since $x_3$ is adjacent to at least 3 vertices in $F$, by the routine argument of ``skew lines" (by \cref{lit}(6)), we can confirm that $x_3$ is only adjacent to vertices in $V_1(F)$. This gives a partition of vertices of $P$. By \cref{decompfunc-linesegment}, $P$ has a line segment for a summand. 
		
		Case 3: $F$ has $2d-1=2(d-1)+1$ vertices. By \cref{dim5:Delta22_2dplus3}, $F$ cannot be $\Delta(2,2)$. By \cite[Theorem 5.9]{PUY2018Excess}, $F$ is either a capped prism or a pentasm. Let $V_1(F), V_2(F)$ be a partition of vertices of $F$ such that the edges between these sets are the projectively parallel edges. Consider the other facet $F'$ that intersects $F$ at a ``side" ridge of $P$. Then $F'$ is decomposable (routine). It follows that $F'$ has a line segment for a summand as well. So there is an edge $[x_1, x_2]$ outside $F$ that is projectively parallel to the ``side" edges of $F$. We may assume without loss of generality that $x_1$ is adjacent to vertices in $V_1(F)$ and $x_2$ is adjacent to vertices in $V_2(F)$. The subgraphs $G[V_1(F)], G[V_2(F)]$ are indecomposable. Thus, the subgraphs $G[x_1 \cup V_1(F)]$, $G[x_2 \cup V_2(F)]$ are indecomposable. It follows that $x_1$ is only adjacent to vertices in $V_1(F)$ and $x_2$ is only adjacent to vertices in $V_2(F)$.  
		
		Now consider any other vertex $x_3$ outside $F$.
		
		Claim 1: If $x_3$ is adjacent to $x_1$, then $x_3$ is not adjacent to $x_2$. 
		
		Suppose on the contrary that $x_3$ is adjacent to both $x_1$ and $x_2$. Since $x_3$ has degree at least 5, it is adjacent to at least 2 vertices in $F$. If $x_3$ has a neighbor in $V_1(F)$, the subgraph $G[x_1 \cup x_3 \cup V_1(F)]$ is indecomposable. Recall that the subgraph $G[x_2 \cup V_2(F)]$ is also indecomposable, so the line containing $[x_2, x_3]$ is skew to some of the ``side" edges of $F$, and connects two indecomposable subgraphs. It follows that there would an indecomposable subgraph that leaves out 1 vertex, making $P$ indecomposable, a contradiction. Thus $x_3$ cannot have a neighbor in $V_1(F)$ and cannot have a neighbor in $V_2(F)$, a contradiction. 
		
		Claim 2: If $x_3$ is adjacent to $x_1$ (or $x_2$, resp.), then the neighbors of $x_3$ in $F$ are only vertices in $V_1(F)$ (or $V_2(F)$, resp.). 
		
		Suppose $x_3$ is adjacent to $x_1$, $v_1, v_2$, where $v_1 \in V_1(F), v_2\in V_2(F)$. Then the subgraph $G[x_1 \cup x_3 \cup V_1(F)]$ is indecomposable and the line containing $[x_3, v_2]$ is skew to some of the ``side" edges of $F$, and connects two indecomposable subgraphs $G[x_1 \cup x_3 \cup V_1(F)]$ and $G[V_2(F)]$, making the union of the two subgraphs indecomposable that leaves out 2 vertices outside. So $P$ would be indecomposable, a contradiction. 
		
		Suppose now $x_3$ is adjacent to at least 2 vertices in $V_2(F)$. Then the subgraph $G[x_3 \cup V_2(F)]$ is indecomposable. The line containing $[x_1, x_3]$ behaves as a skew line to some of the ``side" edges of $F$, which  gives a contradiction as before. This finishes the proof of Claim 2. 
		
		Recall that there are 6 possible subgraphs induced by the 4 vertices outside $F$; see Figure~\ref{fig:prism}. 
		
		In $G_1$, suppose $x_3$ is adjacent to $x_2$ and $x_4$ is adjacent to $x_1$. By Claim 2, the vertex $x_3$ is adjacent to vertices in $V_2(F)$ but not vertices in $V_1(F)$ and the vertex $x_4$ is adjacent to vertices in $V_1(F)$ but not vertices in $V_2(F)$. This gives a partition $V_1=\{ V_1(F), x_1, x_4 \}$, $V_2= \{  V_2(F), x_2, x_3  \}$ of vertices of $P$, and the subgraphs $G[V_1], G[V_2]$ are indecomposable. By \cref{lit}(6), the edges between $V_1$ and $V_2$ are projectively parallel. 
		
		In $G_2$, suppose $x_3$ is adjacent to $x_2$ and $x_4$ is not adjacent to $x_1$. By Claim 2, the vertex $x_3$ is adjacent to vertices in $V_2(F)$ but not vertices in $V_1(F)$. We show that the neighbors of $x_4$ in $F$ are either only vertices in $V_1(F)$ or only vertices in $V_2(F)$. Suppose the contrary. As $x_4$ has degree 1 in $G_2$, the vertex $x_4$ has at least two neighbors in (say) $V_1(F)$. It follows that $G[V_1(F) \cup x_4]$ is indecomposable. Denote by $v$ the neighbor of $x_4$ in $V_2(F)$. Then the line containing $[x_4, v]$ is skew to some of the lines containing the ``side" edges of $F$, making the union of $G[V_1(F) \cup x_4]$ and $V_2(F)$ indecomposable, which leaves out 3 vertices outside. By \cref{lit}(2), the polytope $P$ would be indecomposable, a contradiction. 
		
		If the neighbors of $x_4$ in $F$ are only vertices in $V_1(F)$, the same argument as of the case $G_1$ applies here. If the neighbors of $x_4$ in $F$ are only vertices in $V_2(F)$, this gives a partition $V_1=\{ V_1(F), x_1 \}$, $V_2= \{  V_2(F), x_2, x_3, x_4  \}$ of vertices of $P$, and the subgraphs $G[V_1], G[V_2]$ are indecomposable. By \cref{lit}(6), the edges between $V_1$ and $V_2$ are projectively parallel. 
		
		In $G_3, G_4$, without loss of generality, we may assume that $x_3, x_4$ are both adjacent to $x_1$. By Claim 2, $x_3, x_4$ are only adjacent to vertices in $V_1(F)$. So, we have a partition of vertices of $P$ based on their neighbors in $F$. The result follows as before. And $G_5, G_6$ are not possible by Claim 1. 
		
		Case 4: $F$ has $2d-2=2(d-1)$ vertices. All the decomposable facets are prisms. By the same approach as in Case 3, let $V_1(F), V_2(F)$ be the partition of vertices of $F$ such that the edges between them are the projectively parallel edges. Consider the other facet $F'$ that intersects $F$ at a ``side" ridge of $P$. Then $F'$ is decomposable and is a prism. Thus, there is an edge $[x_1, x_2]$ outside $F$ that is projectively parallel to the ``side" edges of $F$. Note that Claim 1 and Claim 2 work in this case as well. So, if $x_3$ is adjacent to $x_1$, we would know that the neighbors of $x_3$ in $F$ are only the vertices in $V_1(F)$, and if $x_3$ is adjacent to $x_2$, we would know that the neighbors of $x_3$ in $F$ are only the vertices in $V_2(F)$. Suppose $x_3$ is not adjacent to $x_1$ or $x_2$. Then $x_3$ is adjacent to at least 3 vertices in $F$. It is routine to check that the neighbors of $x_3$ in $F$ are either only vertices in $V_1(F)$ or only vertices in $V_2(F)$ (by the skew line argument). The same holds for $x_4, x_5$. So now we know that for every $x_i, i=1,2,3,4,5$, their neighbors in $F$ are either only vertices in $V_1(F)$ or only vertices in $V_2(F)$. These give a suitable partition of vertices of $P$, based on their neighbors. It follows that the edges between the partition are projectively parallel. 
		
	    By \cref{decompfunc-linesegment}, $P$ has a line segment for a summand, and then we can take a suitable cross section of a polytope which is combinatorially equivalent to $P$, where the cross section is either a $(d-1)$-simplex or a $(d-1)$-polytope with $d+1$ vertices. 		
	\end{proof}
\end{theorem}

The three examples given in \cite[Lemma A2]{PUY2018Excess} show that the previous theorem also fails in dimension 4.

The next result remains open in dimension 4, but is false in dimension 3. The decomposability of 3-polytopes with 9 vertices is studied in \cite{BriPrzYos}.

Again, by \cite[Theorem~4]{WJDYconditional}, the structure of decomposable $d$-polytopes with $2d+3$ vertices gives the combinatorial decomposability.

\begin{corollary}\label{cor:2dplus3decomposable}
	There are no conditionally decomposable $d$-polytopes with $2d+3$ vertices, when $d\ge 5$.
\end{corollary}



\begin{remark}
   This manuscript is part of the first author's PhD thesis. Its mathematical significance is twofold: (1) It gave us the confidence and suggested that we should look at decomposable polytopes with up to $3d-4$ vertices \cite{WJDYconditional}; (2) It covers the ``irregular" cases in lower dimensions in \cite{WJDYconditional}.

\end{remark}

\section{Acknowledgement}
 We would like to thank Guillermo Pineda-Villavicencio for carefully reading the first draft and for providing valuable feedback.

\end{document}